\documentclass[11pt,3p,times,a4paper]{elsarticle}

\usepackage{amssymb}

\usepackage[figuresright]{rotating}
\usepackage{natbib,apalike}
\usepackage{makecell,adjustbox}
\usepackage{mathrsfs,mathtools}
\usepackage{algorithm,algorithmicx,algpseudocode}
\usepackage{subfigure}
\usepackage{multirow}
\usepackage{enumitem}
\usepackage{booktabs}
\usepackage{fancyvrb}
\usepackage{amsmath,amssymb,amsfonts,bm,amsthm}
\usepackage{footnote}
\makesavenoteenv{tabular}
\usepackage{setspace}
\usepackage{float}
\usepackage{xcolor}

\usepackage{natbib}
 \bibpunct[, ]{(}{)}{,}{a}{}{,}%
 \def\bibsep{\smallskipamount}%
 \def\BIBand{and}%

\theoremstyle{plain}
\newtheorem{theorem}{Theorem}
\newtheorem{assumption}{Assumption}
\newtheorem{lemma}{Lemma}
\newtheorem{proposition}{Proposition}
\newtheorem{corollary}{Corollary}
\newtheorem{exmp}{Example}

\theoremstyle{remark}
\newtheorem*{remark}{Remark}
\newtheorem*{note}{Note}

\newcommand{\EE}{\mathbb{E}}

\newcommand{\PP}{\mathbb{P}}

\newcommand{\cx}{\mathcal{X}}

\newcommand{\ca}{\mathcal{A}}
\newcommand{\cy}{\mathcal{Y}}
\newcommand{\cn}{\mathcal{N}}

\newcommand{\ck}{\mathcal{K}}
\newcommand{\cl}{\mathcal{L}}
\newcommand{\cz}{\mathcal{Z}}
\newcommand{\cc}{\mathcal{C}}
\newcommand{\cu}{\mathcal{U}}
\newcommand{\cb}{\mathcal{B}}

\newcommand{\cv}{\mathcal{V}}
\newcommand{\cm}{\mathcal{M}}

\newcommand{\cg}{\mathcal{G}}

\newcommand{\scru}{\mathscr{U}}

\newcommand{\bdnew}[1]{\pmb{#1}}
\newcommand{\mfr}{\mathfrak{r}}
\newcommand{\mfm}{\mathfrak{m}}
\newcommand{\mfc}{\mathfrak{c}}
\newcommand{\mfh}{\mathfrak{h}}

\newcommand{\idxbeta}[1]{_{#1\beta}^n}
\newcommand{\sumbeta}{\sum\nolimits_{k=n}^{m(t_n+t)-1}}
\newcommand{\idxalpha}{_{\alpha}^n}
\newcommand{\sumalpha}{\sum\nolimits_{k=n}^{\mu(\tau_n+t)-1}}

\newcommand{\bc}[1]{\bar{c}_{#1}}

\newcommand{\db}[1]{\delta B_{#1}}
\newcommand{\tdb}[1]{\widetilde{\delta B}_{#1}}

\newcommand{\cbrac}[1]{\left\{#1\right\}}
\newcommand{\rbrac}[1]{\left(#1\right)}
\newcommand{\sbrac}[1]{\left[#1\right]}
\newcommand{\norm}[1]{\|#1\|}
\newcommand{\interior}[1]{\operatorname{int}\rbrac{#1}}
\newcommand{\sfomega}{\mathsf{\Omega}}
\newcommand{\kk}[1]{\operatorname{k}\left({#1}\right)}
\newcommand{\holder}{H\"{o}lder's inequality}

\newcommand{\bigo}[1]{O\left(#1\right)}

\newcommand{\ex}[1]{\EE\left[{#1}\right]}

\newcommand{\pr}[1]{\PP\left\{{#1}\right\}}
\newcommand{\I}[1]{\mathbb{I}\left\{{#1}\right\}}

\newcommand{\K}[1]{\operatorname{K}\left({#1}\right)}
\newcommand{\W}[1]{\operatorname{W}_{\mfr}\left({#1}\right)}
\newcommand{\abs}[1]{\left|{#1}\right|}
\newcommand{\F}{\mathscr{F}}

\newcommand{\limn}{\lim_{n\to\infty}}

\makeatletter

\usepackage[hidelinks,colorlinks=false]{hyperref}

\allowdisplaybreaks[4]
\begin{document}

\begin{frontmatter}




\title{A kernel-based stochastic approximation framework for contextual optimization}


\author[fd]{Hao Cao}\ead{hcao21@m.fudan.edu.cn}
\author[fd]{Jian-Qiang Hu\corref{*}}\ead{hujq@fudan.edu.cn}\cortext[*]{Corresponding author.}
\author[sunysb]{Jiaqiao Hu}\ead{jqhu@ams.sunysb.edu}
\address[fd]{School of Management, Fudan University, Shanghai, 200433, China}
\address[sunysb]{Department of Applied Mathematics and Statistics, State University of New York at Stony Brook, Stony Brook, NY 11794}
\date{November 19, 2024}

\begin{abstract}
We present a kernel-based stochastic approximation (KBSA) framework for solving contextual stochastic optimization problems with differentiable objective functions. The framework only relies on system output estimates and can be applied to address a large class of contextual measures, including conditional expectations, conditional quantiles, CoVaR, and conditional expected shortfalls.
Under appropriate conditions, we show the strong convergence of KBSA and characterize its finite-time performance in terms of bounds on the mean squared errors of the sequences of iterates produced. In addition, we discuss variants of the framework, including a version based on high-order kernels for further enhancing the convergence rate of the method and an extension of KBSA for handling contextual measures involving multiple conditioning events.
Simulation experiments are also carried out to illustrate the framework.
\end{abstract}

\begin{keyword}
Simulation \sep Contextual risk \sep Finite difference \sep Kernel smoothing.


\end{keyword}

\end{frontmatter}


\section{Introduction}\label{sec1}
In the era of big data, decision-makers have access not only to system outcomes but also to additional side information, which can significantly enhance the decision-making process. This has led to the development of \emph{contextual stochastic optimization}, where auxiliary information is used to reduce uncertainty and improve decision quality. For instance, in inventory management, decision-makers can use data such as customer demographics, weather forecasts, and economic indicators to better predict product demand and tailor their ordering decisions accordingly
\citep{ban2019big,zhang2024optimal}; whereas in portfolio optimization, historical stock prices and sentiment data from platforms like Twitter can be incorporated to enhance the prediction of stock movements \citep{xu2018stock}. In these settings, data-driven methods are developed to model the conditional distribution of system outcomes, allowing practitioners to make more informed decisions
through optimizing performance measures like conditional expectations \citep{ban2019big, bertsimas2020predictive, kallus2023stochastic,sadana2024survey, wang2025statistical} or conditional quantiles \citep{qi2022distributionally} that reflect the specific context or scenario at hand.
For instance, in the risk-neutral setting, a contextual decision-making problem can be formulated as $\min_{\theta\in\Theta}\EE[c(\theta,Y(\theta))\vert X(\theta)=\bar{x}]$, where $\theta$ denotes the vector of decision variables with a feasible set $\Theta$, $c$ is a cost function, $Y(\theta)$ is the random system outcome, and $X(\theta)$ represents the available contextual information (covariates), which is correlated with $Y(\theta)$. The expectation $\EE[\cdot\vert X(\theta)=\bar{x}]$ is taken with respect to the conditional distribution of $Y(\theta)$ given $X(\theta)=\bar{x}$, where $\bar{x}$ denotes the observed covariate values, and the conditional distribution
$\PP\{Y(\theta)\le y\vert X(\theta)=\bar x\}$ may itself depend on the decision $\theta$.

In the financial domain, the concept of conditional distribution is also crucial in analyzing systemic risk, which arises from the contagion of distress across financial institutions.
A variety of co-risk measures, such as CoVaR (i.e., the quantile of a portfolio loss conditional on another portfolio being at its quantile) and conditional expected shortfall (CoES),
have been developed to study these risks \citep{Adrian2008, acharya2017measuring}. These measures effectively capture the cross-sectional relationship between the performance of an individual institution and the broader financial system through their conditional distributions. There is a growing body of literature on estimating co-risk measures. Some existing approaches include quantile regression \citep{Adrian2008}, time series models \citep{bianchi2023non}, and copula models \citep{karimalis2018measuring}. Because these techniques are parametric and developed based on structural models with specific distributional assumptions, they are often applied on a case-by-case basis, with their effectiveness being heavily reliant on the accuracy of the assumed models.

Unfortunately, estimating and especially optimizing co-risk measures without model-specific assumptions can be extremely difficult. In particular, when continuous random variables are used as covariates, the conditioning events often occur with probability zero, which cannot be directly observed or simulated. In many cases, this is further complicated by the need to estimate the unknown parameters within these conditioning events, e.g., the quantile values in CoVaR and CoES. In addition, unlike some previous studies \citep{ban2019big, kallus2023stochastic, sadana2024survey} where data are assumed to be generated from a fixed distribution, in optimization problems, the distribution of system outcomes and covariates typically depends on the decisions being made (endogenous).
This decision-dependence adds complexity to both sensitivity analysis and optimization of co-risk measures, especially when the underlying distribution is unknown. Finally, the variability in decision variables over time could also lead to nonstationary system outcomes, causing the estimation and optimization errors to be intricately coupled, which may pose additional challenges in algorithm design and analysis.
For stable systems with static input parameters, \cite{huang2022montecarlo} recently develop a non-parametric Monte-Carlo method for CoVaR estimation using order statistics and batching based on simulation samples; \cite{CHH2023} further introduce a stochastic approximation (SA)-type sequential procedure for estimating both CoVaR and its sensitives under a general black-box setting. However, neither of these studies has tackled the more challenging CoVaR optimization problem.

In this paper, we complement prior work by exploring a broad class of \textit{contextual measures}, which encompasses conditional expectations, conditional quantiles, and many co-risk measures such as CoVaR and CoES as special cases. We propose a unifying algorithmic framework, which can not only be applied for simultaneously estimating multiple contextual measures of interest but also allows their sensitivity analysis and optimization to be jointly conducted in a coherent manner. In particular, to address the aforementioned challenges in co-risk measure optimization, our approach employs a kernel smoothing idea \citep{silverman1978weak} to transform an unobservable conditioning event into an unconditional one that is amenable for simulation and uses an SA-type multi-timescale recursive procedure to handle endogenous uncertainty and eliminate estimation errors. The framework we call kernel-based stochastic approximation (KBSA) relies solely on system output estimates, and the number of output evaluations at each iteration is independent of the dimension of the input decision vector, making it suitable for high-dimensional complex systems that cannot be easily modeled using analytical techniques.
Another salient feature of KBSA is that its implementation only requires \textit{two distinct timescales}, regardless of the number of contextual measures to be estimated, allowing for greater flexibility in step-size selection compared to existing multi-timescale sensitivity analysis/optimization algorithms (e.g., \citealp{Hu2022, Hu2023, CHH2023}).
Under appropriate conditions, we prove the almost sure convergence of the method and characterize its finite-time convergence rate by establishing bounds on the mean squared
errors (MSEs) of the sequences of estimates produced by the framework.
Additionally, we present an algorithm acceleration technique based on high-order kernels \citep{schucany1977improvement,fan1992bias} for further enhancing the convergence speed of the method, under which we show that the convergence rate of KBSA is of order $O(n^{-\mfr/(4\mfr+2)})$ for black-box contextual optimization, while the convergence rates for the contextual measures and gradient estimates are of orders $O(n^{-\mfr/(2\mfr+1)})$ and $O(n^{-\mfr/(3\mfr+2)})$, respectively, where $n$ is the number of algorithm iterations and $\mfr\ge2$ is an integer reflecting the smoothness of the system output distribution. Finally, we consider a variant of the framework and discuss how to extend KBSA to handle contextual measures involving multiple conditioning events.

We remark that early research on two-timescale stochastic approximation (SA) primarily focused on asymptotic convergence analysis \cite[cf., e.g.,][]{borkar1997stochastic, borkar2009stochastic, bhatnagar2001two, bhatnagar2003two, tadic2003asymptotic, yin2005discrete}, while the study of convergence rates, particularly finite-time performance, has emerged only more recently. For linear SA methods, \citet{konda2004convergence} analyzed their rate of convergence, which was further refined in subsequent work by \citet{doan2021finite} and \citet{haque2023tight}. For more general nonlinear SA algorithms, \citet{mokkadem2006convergence} established their asymptotic rates, whereas recent studies by \citet{wu2020finite}, \citet{doan2022nonlinear}, and \citet{doan2025fast} examined their finite-time behavior. Nevertheless, the existing finite-time analyses typically rely on assumptions such as unbiasedness and bounded variance of the estimators—conditions that are not directly satisfied in the Kernel-Based Stochastic Approximation (KBSA) framework proposed in this paper, due to the kernel-induced approximation error and potentially unbounded variance. These challenges introduce additional complexities that preclude the straightforward application of existing theoretical results to our setting.


The main contributions of this paper are as follows: 1) We introduce a non-parametric, unifying framework for the estimation, sensitivity analysis, and optimization of a large class of contextual measures; 2) We establish both the strong convergence and the rate of convergence for the KBSA framework. A high-order kernel acceleration method is also proposed to further enhance its performance; 3) We investigate multivariate conditioning events and demonstrate the flexibility and potential of KBSA for a wider range of applications.

The rest of the paper is organized as follows. In Section~\ref{sec.bbo}, we introduce a general formulation of contextual measures and present the KBSA framework for contextual optimization.
The convergence analysis of the proposed algorithmic framework is presented in Section~\ref{sec.converge}. In Section~\ref{sec.rate}, we study the rate of convergence and introduce an acceleration method based on high-order kernels.
In Section~\ref{sec.extend}, the framework is extended to address other conditioning events. Section \ref{sec.experiment} presents numerical experiments evaluating the performance of KBSA, and Section \ref{sec.conclude} concludes the paper.

\section{Blackbox Contextual Optimization}\label{sec.bbo}
\subsection{Notation}
We begin by presenting a glossary of notational symbols frequently used throughout the paper.
\begin{table}[h]
    \caption{A glossary of notation.\label{tab::notation}}
    \centering\small
    {\begin{tabular}{l|l}
    \toprule
    Symbol & Definition\\
    \hline
        $\theta\in\Theta\subset\Re^d$ & vector of decision variables with the feasible region $\Theta$\\
        $Y(\theta)$ & the system outcome\\
        $X(\theta)$ & the covariate, which is correlated with $Y(\theta)$\\
        $\mathcal{Y}(\theta)$ and $\mathcal{X}(\theta)$ & the respective supports of $Y(\theta)$ and $X(\theta)$\\
        $f(\cdot,\theta)$  & marginal probability density function of $X(\theta)$\\
        $F(\cdot,\theta)$ & marginal cumulative distribution function of $X(\theta)$\\
        $\operatorname{K}:\Re\to\Re^+$ & the kernel function, e.g., a probability density\\
        $\lambda(\theta)=(\lambda_1(\theta),\ldots,\lambda_p(\theta))^\intercal$ & vector of contextual measures\\
        $\nu(\theta)$ & the auxiliary measure\\
        $\nabla$ &  gradient operator\\
        $m_j:\Re^{d+j+1}\to\Re$  & user-specified link function, for $j=1,\ldots,p$.\\
        $\bar{m}_j:\Re^{d+j+1}\to\Re$ & $\bar{m}_j(\theta,\lambda_1,\ldots,\lambda_j,\nu):=\EE[m_j(\theta,Y(\theta),\lambda_1,\ldots,\lambda_j)
        \vert X(\theta)\allowbreak=\nu]$, for $j=1,\ldots,p$\\
        $q:\Re^{d+2}\to\Re$ & user-specified link function\\
        $\bar{q}:\Re^{d+1}\to\Re$ &
        $\bar{q}(\theta,\nu):=\EE[q(\theta,X(\theta),\nu)]$\\
        $g:\Re^{d+p}\to\Re$ & user-specified cost function\\
        $a(n) = O(b(n))$ & $\limsup_{n\to \infty}a(n)/b(n)<\infty$ for two positive functions $a(n)$ and $b(n)$\\
        $a(n) = o(b(n))$ & $\limsup_{n\to \infty}a(n)/b(n)=0$ for two positive functions $a(n)$ and $b(n)$\\
        $a(n)=\sfomega(b(n))$ & $\liminf_{n\to\infty}a(n)/b(n)>0$ for two positive functions $a(n)$ and $b(n)$\\
    \bottomrule
    \end{tabular}}
\end{table}

\subsection{Problem Setting}
Let $X(\theta)$ and $Y(\theta)$ be two continuous output random variables from a stochastic system, where $\theta\in \Re^d$ is a parameter vector, whose components may affect the system performance either directly or indirectly through the input distributions. We refer to $X(\theta)$ as the covariate (a.k.a. the feature, attribute, or explanatory variable), which is correlated with the dependent variable $Y(\theta)$. We assume that $(X(\theta), Y(\theta))$ has an absolutely continuous cumulative distribution function (CDF) with a probability density function (PDF). Let $\cx(\theta)$ and $\cy(\theta)$ be the respective supports of $X(\theta)$ and $Y(\theta)$, and  denote by $f(\cdot,\theta)$ and $F(\cdot,\theta)$ the marginal PDF and CDF of $X(\theta)$.

For a given parameter vector $\theta$, we consider the class of \textit{contextual measures} $\lambda(\theta):=(\lambda_1(\theta),\ldots,\lambda_p(\theta))^\intercal\in\Re^p$ that can be expressed as the unique root $\lambda=\lambda(\theta)$ to the following system of nonlinear implicit equations:
\begin{align}
    \EE[m_j(\theta,Y(\theta),
    \lambda_1,\ldots,\lambda_j)\vert X(\theta)=\nu(\theta)]=&0,
    ~~~~\mbox{for $j=1,\ldots,p$,}\label{def.lam}
\end{align}
where $m_j:\Re^{d+j+1}\to\Re$, $j=1,\ldots,p$, are user-defined link functions, and $\nu(\theta)$ denotes an auxiliary measure that is the unique root $\nu=\nu(\theta)$ to another implicit equation:
\begin{align}
    \ex{q(\theta,X(\theta),\nu)}=&0,\label{def.nu}
\end{align}
where $q:\Re^{d+2}\to\Re$ is a user-defined link function.
    Equations~\eqref{def.lam}-\eqref{def.nu} together define an auxiliary problem for characterizing the contextual measures $\lambda(\theta)$.
    For notational convenience, let $\bar{m}_j(\theta,\lambda_1,\ldots,\lambda_j,\nu):=\EE[m_j(\theta,Y(\theta),\lambda_1,\ldots,\lambda_j)
\vert X(\theta)\allowbreak=\nu]$, for $j=1,\ldots,p$, and $\bar{q}(\theta,\nu):=\EE[q(\theta,X(\theta),\nu)]$.
We consider the following two illustrative examples:
\begin{exmp}[Risk-neutral decision-making]\label{exmp_expect}
    Let $c:\Re^{d+1}\to\Re$ be a cost function.
    Define the link functions $m(\theta,y,\lambda):=c(\theta,y)-\lambda$ and $q(\theta,x,\nu):=\bar{x}-\nu$ for some fixed $\bar{x}\in\mathcal{X}$. In this case,  Equations~\eqref{def.lam}-\eqref{def.nu} reduce to
    $$
    \left\{\begin{array}{r}
        \bar{m}(\theta,\lambda,\nu)
        =\mathbb{E}[m(\theta,Y(\theta),
        \lambda)\vert X(\theta)=\nu]=
        \mathbb{E}[c(\theta,Y(\theta))\vert X(\theta)=\nu]
        -\lambda=0,\\
        \bar{q}(\theta,\nu)
        =\mathbb{E}[q(\theta,X(\theta),\nu)]=
        \bar{x}-\nu=0,
    \end{array}\right.
    $$
    whose solution is given by the conditional expected cost  $\lambda=\mathbb{E}[c(\theta,Y(\theta))\vert X(\theta)=\bar{x}]$ and $\nu=\bar{x}$.
\end{exmp}

\begin{exmp}[Systemic risk analysis]\label{exmp_covar}
    Define the link functions
    $m(\theta,y,\lambda):=\psi-\mathbb{I}\{y\le\lambda\}$ and $q(\theta,x,\nu):=\phi-\mathbb{I}\{x\le\nu\}$ for fixed $\psi,\,\phi\in (0,1)$. In this case, Equations~\eqref{def.lam}-\eqref{def.nu} become
    $$
    \left\{\begin{array}{r}
    \bar{m}(\theta,\lambda,\nu)
    =\mathbb{E}[m(\theta,Y(\theta),
    \lambda)\vert X(\theta)=\nu]=
    \psi-\mathbb{P}\{Y(\theta)\le\lambda\vert X(\theta)=\nu\}=0,\\
    \bar{q}(\theta,\nu)
    =\mathbb{E}[q(\theta,X(\theta),\nu)]=
    \phi-\mathbb{P}\{X(\theta)\le \nu \}=0,
    \end{array}\right.
    $$
    whose solution is $\lambda=\inf\{y:\mathbb{P}\{Y(\theta)\le y\vert X(\theta)=\nu\}\ge \psi\}$ and $\nu=\inf\{x:\mathbb{P}\{X(\theta)\le x\}\ge \phi\}$. In this setting, $\lambda(\theta)$ and $\nu(\theta)$ correspond to the CoVaR$_{\phi,\psi}$ of $(X(\theta),Y(\theta))$ and the value-at-risk VaR$_{\phi}$ of $X(\theta)$, respectively.
\end{exmp}

Many existing performance measures such as those mentioned in Section~\ref{sec1} can all be recovered as special cases of $\lambda(\theta)$ through defining appropriate $m_j$ and $q$; see Table~\ref{tab::eg}, which shows examples of some of the contextual measures adopted in the literature and their corresponding link functions. To fix ideas, we focus on the case of a single conditioning event $\{X(\theta)=\nu(\theta)\}$. Variants of the formulation (\ref{def.lam})-(\ref{def.nu}) that involve multiple conditioning events will be discussed in Section~\ref{sec.extend}.
\begin{table}[h]
    \caption{Examples of contextual measures.\label{tab::eg}}
    \centering\small
    {\begin{tabular}{l|l}
    \toprule
    Measures & Link Functions\\
    \hline
        \makecell*[l]{$\lambda_1(\theta)$: the CoVaR$_{\phi,\psi}$ of $(X(\theta),Y(\theta))$;\\
        $\lambda_2(\theta)$: the CoES$_{\phi,\psi}$ of $(X(\theta),Y(\theta))$;\\
        $\nu(\theta)$: the VaR$_\phi$ of $X(\theta)$ \\
        (see \citealp{Adrian2008}).}
        &   \makecell*[l]{
        $m_1(\theta,y,\lambda_1):=\psi-\I{y\le\lambda_1}$;\\
        $m_2(\theta,y,\lambda_1,\lambda_2):=(y-\lambda_2)\I{y\ge\lambda_1}$;\\
        $q(\theta,x,\nu):=\phi-\I{x\le\nu}$.}\\
    \hline
        \makecell*[l]{$\lambda(\theta)=\ex{c(\theta,Y)\vert X=\bar{x}}$: the conditional expected cost \\
        (see \citealp{kallus2023stochastic}).}
        & \makecell*[l]{$m(\theta,y,\lambda):=c(\theta,y)-\lambda$;\\ $q(\theta,x,\nu):=\bar{x}-\nu$, for some given $\bar{x}\in\cx(\theta)$,\\
        where $c:\Re^{d+1}\to\Re$ denote the general cost function.}\\
    \hline
        \makecell*[l]{$\lambda_1(\theta)$: the conditional expectation of $Y(\theta)\vert X(\theta)=\bar{x}$;\\
        $\lambda_2(\theta)$: the conditional variance of $Y(\theta)\vert X(\theta)=\bar{x}$;\\
        $\lambda_3(\theta)$: the $\psi$-conditional quantile of $Y(\theta)\vert X(\theta)=\bar{x}$ \\
        (see \citealp{fan1998efficient, fan2016direct}).}
        &   \makecell*[l]{
        $m_1(\theta,y,\lambda_1):=y-\lambda_1$;\\
        $m_2(\theta,y,\lambda_1,\lambda_2):=(y-\lambda_1)^2-\lambda_2$;\\
        $m_3(\theta,y,\lambda_1,\lambda_2,\lambda_3):=
        \psi-\I{y\le\lambda_3}$;\\
        $q(\theta,x,\nu)=\bar{x}-\nu$, for some given $\bar{x}\in\cx(\theta)$.
        }\\
    \bottomrule
    \end{tabular}
    \begin{note}
        {In the definition of CoVaR$_{\phi,\psi}$ and CoES$_{\phi,\psi}$, $\I{\cdot}$ is the indicator function, and $\phi,\psi\in(0,1)$ are given constants.}
    \end{note}}
\end{table}

Based on the above definition of contextual measures $\lambda(\theta)$, our goal is to solve the following optimization problem:
\begin{align}\label{OPT}
    \min_{\theta\in\Theta}\big\{\bar{g}(\theta):=
    g(\theta,\lambda(\theta))\big\},
\end{align}
where $\Theta\subset \Re^d$ is the feasible region, and
$g:\Re^{d+p}\to\Re$ is a user-specified cost function that may depend on multiple contextual measures. For technical reasons, we assume that $\Theta$ is a convex, compact set specified in terms of inequality constraints $\boldsymbol{c}_j(\theta)\leq 0,~j=1,\ldots,m,$ where each $\boldsymbol{c}_j(\cdot)$ is a continuously differentiable function satisfying $\nabla_{\theta}\boldsymbol{c}_j(\theta)\neq 0$ whenever $\boldsymbol{c}_j(\theta)=0$ \citep[cf., e.g.,][]{kushner2003sa}.

We consider the setting where both the cost function $g$ and the contextual measures $\lambda_1,\ldots,\lambda_p$ are sufficiently smooth, so that a (local) optimal solution to (\ref{OPT}) can be found by employing a gradient descent method of the form
\begin{align}\label{eqn:gd}
    \theta_n=\Pi_\Theta\Big[\theta_{n-1}-\alpha_n\hat{\nabla}_\theta\bar{g}(\theta_{n-1})\Big],~~~~n=1,2,3,\ldots,
\end{align}
where $\alpha_n$ is the step size used at the $n$th iteration, $\Pi_\Theta[\cdot]$ stands for a projection operator that projects an iterate back onto the feasible region $\Theta$ whenever it becomes infeasible,
and $\hat{\nabla}_\theta\bar{g}$ stands for an estimator of the objective function gradient
\begin{align}\label{Eq.gradient}
    \nabla_\theta\bar{g}(\theta)
    =\sum\nolimits_{j=1}^p \frac{\partial g(\theta,\lambda_1,\ldots,\lambda_p)}{\partial \lambda_j}\biggm|_{\lambda=\lambda(\theta)}
    \nabla_\theta \lambda_j(\theta)
    +\nabla_\theta g(\theta,\lambda)|_{\lambda=\lambda(\theta)},~~~~\theta\in\Theta.
\end{align}
Note that central to our setting is that
the distributions defining the expectations in (\ref{def.lam}) and (\ref{def.nu}) are unknown, which precludes an analytical, close-form solution to these equations. Thus, a primary issue is how to reliably estimate the contextual measures $\lambda_j(\theta)$'s and their gradients $\nabla_\theta\lambda_j(\theta)$'s in \eqref{Eq.gradient} based only on system output samples.

\subsection{A Kernel-based Stochastic Approximation Method}\label{sec22}
As mentioned in Section~\ref{sec1}, one challenge in directly approximating the expectation in (\ref{def.lam}) is that it requires data on $Y(\theta)$ to be collected from a conditional distribution given that the covariate $X(\theta)$ is held constant at $\nu(\theta)$, a zero-probability scenario that is impossible to be replicated or even simulated. This is in addition to the fact that the exact value of $\nu(\theta)$ is typically unknown, which itself needs to be estimated through the implicit equation (\ref{def.nu}).
Our main idea to overcome these difficulties is to use a kernel-smoothing technique to transform the conditional expectation in (\ref{def.lam}) into an unconditional one that can be easily simulated. In particular, under mild technical conditions, the following relationship can be derived from Bochner’s lemma \citep{parzen1962estimation}:
\begin{align}\label{eqn:tr}
    \lim_{h\to0^+}
    \ex{\frac{1}{h}\K{\frac{\nu-X}{h}}
    m_j(\theta,Y,
    \lambda_1,\ldots,\lambda_j)}
    =f(\nu,\theta)\ex{m_j(\theta,Y,\lambda_1,\ldots,\lambda_j)\vert X=\nu},
\end{align}
for $j=1,\ldots,p$, where $\K{\cdot}$ standards for a kernel function (e.g., a probability density) and
$h>0$ is the associated kernel-bandwidth parameter.
Assuming that $f(\nu,\theta)>0$, this suggests that solving Equation~\eqref{def.lam} is  equivalent to finding the solution to the stochastic root-finding problem
\begin{align}\label{Eq.kernel-smoothing}
    \lim_{h\to0^+}
    \ex{\frac{1}{h}\K{\frac{\nu(\theta)-X(\theta)}{h}}
    m_j(\theta,Y(\theta),
    \lambda_1(\theta),\ldots,\lambda_j(\theta))}
    =0~~~~\mbox{for $j=1,\ldots,p$}.
\end{align}
As contrasted with \eqref{def.lam}, the expectation in (\ref{Eq.kernel-smoothing}) can be conveniently estimated by sampling from the joint distribution of $(X(\theta),Y(\theta))$,
while the limit can be approximated by employing a suitable sequence of diminishing bandwidths $\{h_n\}$.
Consequently, given the values of $\nu(\theta)$, $\lambda_1(\theta),\ldots,\lambda_{j-1}(\theta)$, the $j$th contextual measure $\lambda_j(\theta)$ can be recursively estimated by applying the standard root-finding SA algorithm for approximately solving the $j$th equation in (\ref{Eq.kernel-smoothing}), i.e.,
\begin{align}\label{eqn:sa1}
    \lambda_{j,n}=
    \lambda_{j,n-1}+\frac{\beta_{n}}{h_n}\K{\frac{\nu(\theta)-X(\theta)}{h_n}}
    m_j(\theta,Y(\theta),\lambda_{1}(\theta),\ldots,\lambda_{j-1}(\theta),\lambda_{j,n-1}),
\end{align}
for $n=1,2,\ldots$, where $(X(\theta),Y(\theta))$ is a pair of system output random variables generated under the input parameter vector $\theta$, $\lambda_{j,n}$ denotes an estimate of $\lambda_j$ obtained at iteration $n$, and $\beta_{n}$ and $h_n$ are the respective step size and bandwidth used.
Note that the same recursion can be carried out simultaneously to jointly estimate $\nu(\theta)$ and other contextual measures in $\lambda(\theta)$, with $\nu(\theta)$ and $\lambda_1(\theta),\ldots,\lambda_{j-1}(\theta)$ on the right-hand side of (\ref{eqn:sa1}) replaced by their current estimates $\nu_{n-1}$, $\lambda_{1,n-1},\ldots,\lambda_{j-1,n-1}$.
This then yields a high-dimensional version of (\ref{eqn:sa1}) (see recursion~\eqref{Eq.nu} in the proposed KBSA algorithm), which can be viewed as an SA method for concurrently solving a system of $p+1$ zero-finding problems.

A similar but slightly more sophisticated procedure can also be applied to estimate the contextual measure gradients. For notational convenience, we have defined
$
    \bar{m}_j(\theta,\lambda_1,\ldots,\lambda_j,\nu)
    :=\EE[m_j(\theta,Y,\lambda_1,\ldots,\lambda_j)\vert X=\nu],
$
and can write Equation~\eqref{def.lam} as $\bar{m}_j(\theta,\lambda_1(\theta),\ldots,\lambda_j(\theta),\nu(\theta))=0$ for $j=1,\ldots,p$. Thus, if we assume the smoothness of $\bar{m}_j$, then by directly differentiating both sides of this equation using the chain rule, we arrive at the following implicit equation characterizing  $\nabla_{\theta}\nu(\theta)$ and $\nabla_\theta \lambda_j(\theta)$'s:
\begin{align}\label{def.gradients}
    \nabla_\theta\bar{m}_j(\theta,\lambda_1,\ldots,\lambda_j,\nu)
    |_{\lambda_1=\lambda_1(\theta),\ldots,
    \lambda_j=\lambda_j(\theta),\nu=\nu(\theta)}
    +\sum\nolimits_{i=1}^j\frac{\partial\bar{m}_j}{\partial
    \lambda_i}\nabla_\theta\lambda_i(\theta)
    +\frac{\partial\bar{m}_j}{\partial \nu}
    \nabla_\theta\nu(\theta)=0.
\end{align}
Unfortunately, the equation cannot be readily solved, because the calculation of the partial derivatives of $\bar{m}_j$ requires explicit knowledge of the system output distribution, which is not available in our setting. So we instead propose an approach based on finite-difference (FD) approximation of the entire left-hand side of \eqref{def.gradients}.
To illustrate the main idea underlying the approach, we begin with the case where $\theta$ is a scalar and consider the conditional expectation function $\bar m_j$ evaluated at its slightly perturbed input values: $\theta^\pm:=\theta\pm c$, $\lambda_i^\pm:=\lambda_i(\theta)\pm c\lambda_i'(\theta)$ for $i=1,\ldots,j$, and $\nu^\pm=\nu(\theta)\pm c\nu'(\theta)$. We define $\bar{m}_j^\pm:=\bar{m}_j(\theta^\pm,\lambda_1^\pm,\ldots,\lambda_j^\pm,\nu^\pm)$ and $m_j^\pm:=m_j(\theta^\pm,Y^\pm,
\lambda_1^\pm,\ldots,\lambda_j^\pm)$. In particular, through a straightforward Taylor series expansion of $\bar m_j$ around $\theta$, $\lambda_i(\theta)$ and $\nu(\theta)$, we have that
\begin{align*}
    \lim_{c\to0}\left(
    \frac{\bar{m}_j^+-\bar{m}_j^-}{2c}
    \right)
    =\frac{\partial\bar{m}_j}{\partial \theta}
    +\sum\nolimits_{i=1}^j\frac{\partial\bar{m}_j}{\partial \lambda_i}\lambda_i'(\theta)
    +\frac{\partial\bar{m}_j}{\partial \nu}\nu'(\theta).
\end{align*}
Next, applying the kernel smoothing equation (\ref{eqn:tr}) to each of the perturbed $\bar m_j$ functions in the difference quotient, it can be derived that
\begin{align}
    \lim_{c\to0}\lim_{h\to0^+}
    \ex{\frac{1}{h^2}\K{\frac{\nu^+-X^+}{h}}
    \K{\frac{\nu^--X^-}{h}}
    \frac{m_j^+-m_j^-}{2c}
    }
    =&\lim_{c\to0}f(\nu^+,\theta^+)f(\nu^-,\theta^-)\left(
    \frac{\bar{m}_j^+-\bar{m}_j^-}{2c}
    \right)\nonumber\\
    =&f^2(\nu,\theta)\left(
    \frac{\partial\bar{m}_j}{\partial \theta}
    +\sum\nolimits_{i=1}^j\frac{\partial\bar{m}_j}{\partial \lambda_i}\lambda_i'(\theta)
    +\frac{\partial\bar{m}_j}{\partial \nu}\nu'(\theta)\right),\label{KBFD}
\end{align}
where $(X^\pm,Y^\pm):=(X(\theta^\pm),Y(\theta^\pm))$ are two pairs of output random variables independently generated under the perturbed parameter vectors $\theta^\pm$. Notice that
the expression in the parenthesis of (\ref{KBFD}) is precisely the left-hand side of (\ref{def.gradients}) in one dimension. This suggests that once the values of $\nu(\theta)$, $\lambda_1(\theta),\ldots,\lambda_j(\theta)$ are determined, their derivatives
$\nu'(\theta)$ and $\lambda_1'(\theta),\ldots,\lambda_{j}'(\theta)$ can be jointly obtained as the solution to the following system of stochastic root-finding equations:
{\small
\begin{align}\label{eqn:tmp}
    \lim_{c\to0}\lim_{h\to0^+}
    \ex{\frac{1}{h^2}\K{\frac{\nu^+-X^+}{h}}
    \K{\frac{\nu^--X^-}{h}}
    \frac{m_j^+-m_j^-}{2c}
    }=0,
\end{align}}
for $j=1,\ldots,p.$ These equations can then be approximately solved by using SA-type of recursions similar to (\ref{eqn:sa1}) in a way that is analogous to solving (\ref{Eq.kernel-smoothing}).

Because the usual FD perturbation is carried out element-wise, in the multi-dimensional ($d\ge1$) case, each iteration of the resultant SA recursions for solving (\ref{eqn:tmp}) would require
$2d$ output pair evaluations.
This could be computationally demanding for high-dimensional problems, especially when the output evaluations are expensive. Following \cite{spall1992multivariate}, we therefore consider an alternative simultaneous perturbation (SP)-based implementation, where the idea is to replace the $d$ unit directions in the FD method with a single random direction vector $\Delta_n\in\Re^d$ that simultaneously varies all components of an input vector in random directions.
Specifically, let $\nu_n$, $\lambda_{j,n}$, $G_{\nu,n}$, $G_{j,n}$ be the respective estimates of $\nu(\theta_n)$, $\lambda_j(\theta_n)$, $\nabla_\theta\nu(\theta_n)$, $\nabla_\theta\lambda_j(\theta_n)$, and to simplify notation, define
{\begin{align}
    &\bc{n}=c_n/\max\cbrac{1,
    \norm{G_{\nu,n-1}}/\sqrt{d},
    \norm{G_{1,n-1}}/\sqrt{d},\ldots,
    \norm{G_{p,n-1}}/\sqrt{d}
    },\nonumber\\
    &\theta_n^\pm:=\theta_n\pm \bc{n+1}\Delta_{n+1},~~~~~~~
    (X_n^\pm,Y_n^\pm):=(X(\theta_n^\pm),Y(\theta_n^\pm)),\nonumber\\
    &\nu_n^\pm:=\nu_n\pm \bc{n+1}\Delta_{n+1}^\intercal G_{\nu,n},~~
    \lambda_{j,n}^\pm:=\lambda_{j,n}\pm \bc{n+1}\Delta_{n+1}^\intercal G_{j,n},\label{tmp2}\\
    &K_n:=\frac{1}{h_{n+1}}\K{\frac{\nu_n-X(\theta_n)}{h_{n+1}}},~~
    K_n^\pm:=\frac{1}{h_{n+1}}\K{\frac{\nu_n^\pm-X_n^\pm}{h_{n+1}}},
    \label{term_Kernel_uni}\\
    &\bdnew{\lambda}_n:=(\lambda_{1,n},\lambda_{2,n},\ldots,\lambda_{p,n})^\intercal,~~
    \bdnew{\lambda}_n^\pm:=(\lambda_{1,n}^\pm,\lambda_{2,n}^\pm,\ldots,\lambda_{p,n}^\pm)^\intercal.\nonumber
\end{align}}
The SP version of the SA recursion for solving (\ref{eqn:tmp}) in high dimensions is presented in (\ref{Eq.Gnu}) of Algorithm 1, where to prevent the perturbations in $\nu_n^\pm$ and $\lambda_{j,n}^\pm$ from becoming excessive large (see (\ref{tmp2})), we have used a
perturbation size $\bar c_n$ that scales $c_n$ by the magnitudes of the gradient estimates \citep[cf., e.g.,][]{Hu2023}. Note that as compared with the conventional FD implementation, each iteration of (\ref{Eq.Gnu}) only requires two output pair samples, independent of the problem dimension.

\begin{algorithm*}
\caption{KBSA for Blackbox Contextual Optimization.\label{alg}}
\begin{algorithmic}
   \State\textbf{Input} initial estimates $\nu_0,\bdnew{\lambda}_0,G_{\nu,0},G_{1,0},\ldots,G_{p,0},\theta_0$; step-size sequences $\{\alpha_n\}$ and $\{\beta_n\}$; a perturbation-size sequence $\{c_n\}$; a kernel-bandwidth sequence $\{h_n\}$;
   \State  \textbf{Initialize} the iteration counter $n\gets 1$;
   \State  \textbf{Iterate} until a stopping rule is satisfied:
   {\allowdisplaybreaks\footnotesize
   \begin{flalign}
   \hspace{-4mm}
    \begin{pmatrix}
        \nu_n\\\bdnew{\lambda}_n
    \end{pmatrix}=&
    \begin{pmatrix}
        \nu_{n-1}\\\bdnew{\lambda}_{n-1}
    \end{pmatrix}
    +\beta_n
    \begin{pmatrix*}[l]
        q(\theta_{n-1},X(\theta_{n-1}),\nu_{n-1})\\
        m_1(\theta_{n-1},Y(\theta_{n-1}),\lambda_{1,n-1})K_{n-1}\\ 
        m_2(\theta_{n-1},Y(\theta_{n-1}),\lambda_{1,n-1},\lambda_{2,n-1})K_{n-1}\\ \vdots\\
        m_p(\theta_{n-1},Y(\theta_{n-1}),\bdnew{\lambda}_{n-1})K_{n-1}
    \end{pmatrix*},\label{Eq.nu}&\\
    \hspace{-4mm}
    \begin{pmatrix}
        G_{\nu,n}\\G_{1,n}\\ \vdots\\G_{p,n}
    \end{pmatrix}=&
    \begin{pmatrix}
        G_{\nu,n-1}\\G_{1,n-1}\\ \vdots\\G_{p,n-1}
    \end{pmatrix}
    +\beta_n
    \begin{pmatrix*}[l]
    \big[q(\theta_{n-1}^+,X_{n-1}^+,\nu_{n-1}^+)
    -q(\theta_{n-1}^-,X_{n-1}^-,\nu_{n-1}^-)\big]/(2\bc{n}\Delta_n)\\
    \big[m_1(\theta_{n-1}^+,Y_{n-1}^+,\lambda_{1,n-1}^+)
    -m_1(\theta_{n-1}^-,Y_{n-1}^-,\lambda_{1,n-1}^-)\big]
    K_{n-1}^+K_{n-1}^-/(2\bc{n}\Delta_n)\\ 
    \big[m_2(\theta_{n-1}^+,Y_{n-1}^+,\lambda_{1,n-1}^+,\lambda_{2,n-1}^+)
    -m_2(\theta_{n-1}^-,Y_{n-1}^-,\lambda_{1,n-1}^-,\lambda_{2,n-1}^-)\big]
    K_{n-1}^+K_{n-1}^-/(2\bc{n}\Delta_n)\\ \vdots\\
    \big[m_p(\theta_{n-1}^+,Y_{n-1}^+,\bdnew{\lambda}_{n-1}^+)
    -m_p(\theta_{n-1}^-,Y_{n-1}^-,\bdnew{\lambda}_{n-1}^-)\big]
    K_{n-1}^+K_{n-1}^-/(2\bc{n}\Delta_n)\\ 
    \end{pmatrix*},\label{Eq.Gnu}&\\
    \theta_n=&\Pi_\Theta\Big[\theta_{n-1}-\alpha_n
        \Big(\nabla_\theta g(\theta_{n-1},\bdnew{\lambda}_{n-1})
        +\sum\nolimits_{j=1}^p\frac{\partial g(\theta_{n-1},\lambda_{1,n-1},\ldots,\lambda_{p,n-1})}{\partial \lambda_j}G_{j,n-1}\Big)\Big],\label{Eq.theta}&\\
        n\gets& n+1.\nonumber&
   \end{flalign}}
\end{algorithmic}
\end{algorithm*}

Our proposed KBSA framework incorporates the contextual measure and their gradient estimates generated by (\ref{Eq.nu}) and (\ref{Eq.Gnu}) into the gradient descent method (\ref{eqn:gd}) in a straightforward way and is presented in Algorithm 1.

It is not difficult to observe that because the constructions of the perturbed values  $\nu_{n-1}^\pm$ and $\lambda_{i,n-1}^\pm$ rely on the estimates $\nu_{n-1}$ and $\lambda_{j,n-1}$ (see, e.g., (\ref{tmp2})),
the two processes \eqref{Eq.nu} and \eqref{Eq.Gnu} for estimating contextual measures and their gradients are intrinsically coupled. Unlike existing studies, where this coupling issue is addressed through the choice of distinct step-sizes for different recursions  \citep[e.g.,][]{Hu2022,CHH2023,Hu2023}, our approach instead views (\ref{Eq.kernel-smoothing}) and (\ref{eqn:tmp}) jointly as a single root-finding problem, so that
both \eqref{Eq.nu} and \eqref{Eq.Gnu} become parts of the same SA recursion for finding the solution to the problem.
As a result, although KBSA comprises three separate recursions \eqref{Eq.nu}-\eqref{Eq.theta}, it is in effect a two-timescale SA method, regardless of the number of contextual measures involved in the objective function. This feature of the method allows for greater flexibility in step-size selection compared to SA algorithms operating on three or more timescales \citep{borkar2009stochastic}.

In the implementation of KBSA, the step size $\alpha_n$ should be chosen very small relative to $\beta_n$ to ensure that \eqref{Eq.nu} and \eqref{Eq.Gnu} are carried out at a rate that is much faster than \eqref{Eq.theta}. Intuitively, this is because the contextual measures and their gradients are computed with respect to a fixed parameter vector $\theta_{n}$. Thus, using a small (vanishing) $\alpha_n$ has the effect of making the changes in $\theta_{n}$ to become negligible, so that in the long run, the sequence $\{\theta_n\}$ would appear to be fixed at a constant value when viewed from \eqref{Eq.nu} and \eqref{Eq.Gnu}.
We note that for ease of exposition, both \eqref{Eq.nu} and \eqref{Eq.Gnu} are presented using the same $\beta_n$. However, this does not necessarily mean that all components of the iterates in \eqref{Eq.nu} and \eqref{Eq.Gnu} should be updated based on the same step-size.
In practice, due to the difference in the magnitude of various estimates, it could be beneficial to consider using different constant multipliers in the step-size for different components to improve the empirical performance of the method (see, e.g., Section~\ref{sec.experiment}).

\section{Convergence Analysis}\label{sec.converge}
Let $(\Omega, \F, \PP)$ be the probability space induced by the algorithm, where $\Omega$ is the collection of all sample trajectories that could possibly be generated by executing the algorithm, $\F$ is a $\sigma$-field of subsets of $\Omega$, and $\PP$ is a probability measure induced by the algorithm. Let $\F_n:=\sigma\{(\nu_k, \bdnew{\lambda}_{k}, {G}_{\nu,k}, {G}_{1,k},\ldots,G_{p,k},\theta_k):0\leq k\leq n\}$ be the family of increasing $\sigma$-fields generated by the set of iterates obtained up to iteration $n=0,1,\ldots$.
For a given column vector $\lambda\in\Re^p$, $\lambda^\intercal$ denotes the transpose and $\|\lambda\|$ denotes the Euclidean norm.
For two positive functions $a(n)$ and $b(n)$, we write
$a(n) = O(b(n))$ if $\limsup_{n\to \infty}a(n)/b(n)<\infty$,
$a(n) = o(b(n))$ if $\limsup_{n\to \infty}a(n)/b(n)=0$,
and $a(n)=\sfomega(b(n))$ if $\liminf_{n\to\infty}a(n)/b(n)>0$.
For a given set $A$, let $\interior{A}$ denote the interior of $A$.

We begin by stating the assumptions that will be used in our analysis.
\begin{assumption}\label{A.stats}
For almost all $(\theta_n,\nu_n,\bdnew{\lambda}_n)$ pairs, there exists an open neighborhood $\cn\times\cx_\varepsilon\times\Lambda_\varepsilon$ of $(\theta_n,\nu_n,\bdnew{\lambda}_n)$, independent of $n$ and $\omega\in\Omega$, such that
\begin{enumerate}[itemsep=0pt, parsep=0pt,
label=\upshape(\alph*), ref=\theassumption (\alph*)]
    \item\label{A.marginal_distribution}  $f(x,\theta)$ is three-times continuously differentiable
    in $x$ on $\interior{\cx(\theta)}$ for every $\theta\in\cn$ and in $\theta$ on $\cn$, with uniformly bounded derivatives;
    \item\label{A.marginal_distribution2}
    $f(\cdot,\cdot)\ge \varepsilon$ on $\cx_\varepsilon\times\cn$ for some constant $\varepsilon>0$;
    \item\label{A.q}
    $\bar{q}(\theta,\nu):=\ex{q(\theta,X(\theta),\nu)}$ is three-times continuously differentiable
    in $\nu$ on $\interior{\cx(\theta)}$ for every $\theta\in\cn$ and in $\theta$ on $\cn$, with uniformly bounded derivatives;
    \item\label{A.q2}
    $\EE[q^2(\theta,X(\theta),\nu)]\le C_q$ and
    $\partial\bar{q}(\theta,\nu)/\partial \nu\le-\varepsilon_q$ for  $(\theta,\nu)\in\cx_\varepsilon\times\cn$ and some constants $C_q,\varepsilon_q>0$;
    \item\label{A.m} $\bar{m}_j(\theta,\lambda_1,\ldots,\lambda_j,\nu):=\EE[m_j(\theta,Y(\theta),\lambda_1,\ldots,\lambda_j)\vert  X(\theta)\allowbreak=\nu]$ is three-times continuously differentiable
    in $\nu$ on $\interior{\cx(\theta)}$ for every $\theta\in\cn$, and in $(\theta,\lambda)$ on $\cn\times\Lambda_\varepsilon$,
    for $j=1,\ldots,p$, with uniformly bounded derivatives;
    \item\label{A.m2}
    $\EE[m_j^2(\theta,Y(\theta),\lambda_1,\ldots,\lambda_j)\vert X(\theta)=\nu]\le C_m$ and
    $\partial \bar{m}_j(\theta,\lambda_1,\ldots,\lambda_j,\nu)/\partial \lambda_j\le-\varepsilon_m$
    for all $(\theta,\nu,\lambda)\in\cn\times\cx_\varepsilon\times\Lambda_\varepsilon$,
    some constants $C_m,\varepsilon_m>0$, and $j=1,\ldots,p$;
    \item\label{A.X-Lipschitz} there exists an r.v. $\kappa$ with $\ex{\kappa}<\infty$ such that
    $\abs{X(\theta_1)-X(\theta_2)}\leq \kappa\norm{\theta_1-\theta_2}$ almost surely (a.s.) for all $\theta_1,\theta_2\in\cn$;
    \item\label{A.g-conti}
    $g$ is continuously differentiable with uniformly bounded derivatives;
    \item\label{A.Schwartz}
    $f(\cdot,\theta)$ and $\bar{m}_j(\theta,\lambda_1,\ldots,\lambda_j,\cdot)$ are Schwartz functions for every given $(\theta,\lambda)\in\cn\times\Lambda_\varepsilon$ and $j=1,\ldots,p$.
\end{enumerate}
\end{assumption}

\begin{assumption}\label{A.direction}
    The random directions $\{\Delta_{n}\}$ are i.i.d., independent of $\F_{n-1}$. The components of $\Delta_{n}$ are mutually independent and follow the Bernoulli distribution $\mathbb{P}\{\Delta_{n,j}=1\}=\mathbb{P}\{\Delta_{n,j}=-1\}=1/2$ for $j\in\{1,\ldots,d\}$ and $n\geq 1$.
\end{assumption}

\begin{assumption}\label{A.stepsize_uni}
The step-, perturbation-sizes, and bandwidths satisfy the following conditions:
\begin{enumerate}[itemsep=0pt, parsep=0pt,
label=\upshape(\alph*), ref=\theassumption (\alph*)]
    \item\label{A.alpha_uni} $\alpha_n>0$, $\sum_{n=1}^\infty\alpha_n=\infty$;
    \item\label{A.beta_uni} $\beta_n,c_n,h_n>0$, $\limn c_n=0$, $\limn h_n=0$,
    $\sum_{n=1}^\infty \beta_n^2/(c_nh_n)^2<\infty$,
    $\beta_n=o(c_n^2h_n^2)$;
    \item\label{A.multi_timescale_uni} $\alpha_n=o(\beta_n)$.
\end{enumerate}
\end{assumption}

\begin{assumption}\label{A.kernel}~\
\begin{enumerate}[itemsep=0pt, parsep=0pt,
label=\upshape(\alph*), ref=\theassumption (\alph*)]
    \item\label{A.kernel-function} The kernel function $\operatorname{K}$ is a symmetric, bounded, and continuous density function. It is absolutely integrable with respect to the Lebesgue measure and satisfies
    $
        \lim_{\abs{u}\to\infty}u\K{u}=0,~
        \int_{\Re}\abs{u\ln{\abs{u}}}^{1/2}
        \abs{d\K{u}}<\infty.
    $
    \item\label{A.kernel-cont} There exists a constant $k_0>0$ such that $\abs{\K{u_1}-\K{u_2}}\leq k_0\abs{u_1-u_2}$ for all $u_1,u_2\in\Re$.
    \item\label{A.kernel-pos} There exists a constant $\hat{\varepsilon}>0$ and an open neighborhood $\cu$ of $0$ such that $\K{u}\geq\hat{\varepsilon}$ for all $u\in \cu$.
    \item\label{A.kernel-transform}
    The Fourier transform of the kernel function
    $\kk{u}=\int e^{-iuy}\K{y}dy,~u\in\Re,$
    is absolutely integrable and satisfies
    $\lim_{u\to0}[1-\kk{u}]/u^2:= k^*<\infty$.
\end{enumerate}
\end{assumption}

\begin{remark}
    Assumption~\ref{A.stats} states the regularity conditions on the marginal density function $f$, and the link functions $q$ and $m_j$.
    The conditions imposed on the density function $f$ are typically satisfied by many standard continuous distributions encountered in practice—such as the normal, exponential, and uniform distributions—and are commonly assumed in the literature on risk control \citep[cf., e.g.,][]{Hu2022, Hu2023} and kernel estimation \citep[e.g.,][]{silverman1978weak}.
    Note that we only require the continuity and boundedness of the expectations $\bar{q}$ and $\bar{m}_j$ rather than those of the original functions $q$ and $m_j$, which is inline with the conditions used in, e.g.,  \cite{spall1992multivariate, CHH2023, Hu2023}.
    For instance, in the case of risk-neutral contextual optimization (see Example~\ref{exmp_expect}), we have $\bar{m}(\theta,\lambda,\nu)=\EE[c(\theta,Y(\theta))|X(\theta)=\bar{x}]-\lambda$ for some given $\bar{x}$,  
    and the continuity and boundedness of $\bar{m}$ hold under mild regularity conditions on the distribution and the cost function \citep[cf.,][]{wang2025statistical}.
    In the CoVaR setting (see Example~\ref{exmp_covar}), we have $\bar{q}(\theta,\nu)=\phi-\pr{X(\theta)\le \nu}$ and $\bar{m}(\theta,\lambda,\nu)=\psi-\pr{Y(\theta)\le\lambda\vert X(\theta)=\nu}$, both of which are automatically bounded. Moreover, their continuity can be readily ensured under mild assumptions on the underlying distribution \citep[cf.,][]{CHH2023, huang2022montecarlo}.
    Assumption~\ref{A.stats} implies that the functions $\lambda(\theta)$ and $\nu(\theta)$ are locally Lipschitz continuous in $\theta$. For example, under Assumptions~\ref{A.q} and \ref{A.m}, 
    by differentiating both sides of the equations $\bar{m}_1(\theta,\lambda_1(\theta),\nu(\theta))=0$ and $\bar{q}(\theta,\nu(\theta))=0$ with respect to $\theta$ and applying the chain rule, we find that
    \begin{align*}
        \nabla_\theta\bar{m}_1(\theta,\lambda_1,\nu)|_{\lambda_1=\lambda_1(\theta),\nu=\nu(\theta)}
        +\frac{\partial \bar{m}_1(\theta,\lambda_1(\theta),\nu(\theta))}{\partial\lambda_1}\nabla_\theta\lambda_1(\theta)
        +\frac{\partial \bar{m}_1(\theta,\lambda_1(\theta),\nu(\theta))}{\partial\nu}\nabla_\theta\nu(\theta)=0,\\
        \nabla_\theta\bar{q}(\theta,\nu)|_{\nu=\nu(\theta)}
        +\frac{\bar{q}(\theta,\nu(\theta))}{\partial\nu}\nabla_\theta\nu(\theta)=0.
    \end{align*}
    Together with Assumptions~\ref{A.q}-\ref{A.m2}, these expressions indicate the existence of constants $C_{\lambda},C_{\nu}>0$ such that $\|\nabla_\theta\lambda_1(\theta)\|\le C_{\lambda}$ and $\|\nabla_\theta\nu(\theta)\|\le C_{\nu}$ for all $\theta\in\mathcal{N}$. As a result, by the mean value theorem, for any $\theta_1,\theta_2\in\mathcal{N}$, we have
        $|\lambda_1(\theta_1)-\lambda_1(\theta_2)|
        =|(\theta_1-\theta_2)^\intercal
        \nabla_\theta\lambda_1(\theta)|_{\theta=\bar\theta}|
        \le C_\lambda \|\theta_1-\theta_2\|,$
    for some $\bar\theta\in\mathcal{N}$. Analogously, we can show that $\nu(\theta)$ and $\lambda_j(\theta)$ ($j=1,\ldots,p$) are locally Lipschitz continuous in $\theta$.
    Assumptions~\ref{A.direction} and \ref{A.stepsize_uni} are standard in the SA literature (e.g., \citealp{spall1992multivariate, kushner2003sa, hu2024, li2024eliminating, cao2025infinitesimal}). As discussed in Section~\ref{sec22}, Condition~\ref{A.multi_timescale_uni} ensures that the recursions in KBSA are carried out at distinct timescales. This condition, together with \ref{A.alpha_uni}-\ref{A.beta_uni}, implies that $\sum_{n=1}^\infty\alpha_n^2<\infty$ and $\sum_{n=1}^\infty\beta_n=\infty$. 
    Assumption~\ref{A.kernel} is also standard conditions used in the kernel estimation/regression literature to ensure the validity of the kernel approximation technique (cf. e.g., \citealp{parzen1962estimation}).  These conditions are satisfied by many kernel functions, e.g., Gaussian densities of the form $\K{u}=(2\pi)^{-1/2} e^{-u^2/2}$.
\end{remark}

The convergence of KBSA is investigated based on an ordinary differential equation (ODE) method (e.g., \citealp{kushner2003sa, borkar2009stochastic}).
In Section~\ref{sec.preliminary}, we begin by establishing the boundedness of the sequence of iterates $\{(\nu_n,\bdnew{\lambda}_n,G_{\nu,n},G_{1,n},\allowbreak \ldots,G_{p,n},\theta_n)\}$ in Lemmas~\ref{LM_bound_nu}-\ref{LM_bound_G}, which allows us to construct continuous-time interpolations of these iterates.
Then in Sections~\ref{sec.converge-measure} and \ref{sec.converge-theta}, we proceed to analyze the properties of these interpolated iterates and use a system coupled ODEs to characterize their asymptotic behavior. This finally leads us to
conclude that the sequence $\{\theta_n\}$ generated by \eqref{Eq.theta} approaches the limiting solution to the following projected ODE:
\begin{align}\label{ODE_theta}
    \dot{\theta}(t)=-\nabla_\theta\bar{g}(\theta)|_{\theta=\theta(t)}+z(t),~~~~t\ge0,
\end{align}
where $\nabla_\theta\bar{g}$ is given by \eqref{Eq.gradient} and $z(t)\in-\cc(\theta(t))$ is the minimum force needed to keep the trajectory $\theta(t)$ from leaving the constraint set $\Theta$ with $\cc(\theta):=\{x\in\Re^d:x^\intercal(\vartheta-\theta)\le0,~\forall\vartheta\in\Theta\}$ being the normal cone of $\Theta$ at $\theta$.

\subsection{Preliminary Results}\label{sec.preliminary}
We first present Lemmas~\ref{LM_bound_nu}-\ref{LM_bound_G}, which show that the sequence of estimates $\{(\nu_n,\bdnew{\lambda}_n,G_{\nu,n},G_{1,n},\ldots,G_{p,n})\}$ generated by \eqref{Eq.nu}-\eqref{Eq.Gnu}
remains bounded at all times, both in the second moment and almost surely (a.s.). Their proofs are given in supplementary materials.

\begin{lemma}\label{LM_bound_nu}
    If Assumptions~\ref{A.q}-\ref{A.q2},
    \ref{A.beta_uni}, and \ref{A.multi_timescale_uni} hold, then
        \upshape{(a)} $\sup_n\mathbb{E}[\nu_n^2]<\infty$;~~(b) $\sup_n|\nu_n|<\infty$ a.s.
\end{lemma}

\begin{lemma}\label{LM_bound_lam}
    If Assumptions~\ref{A.marginal_distribution}, \ref{A.marginal_distribution2}, \ref{A.m}, \ref{A.m2}, \ref{A.Schwartz}, \ref{A.beta_uni}
    and \ref{A.kernel} hold, then
    \upshape{(a)} $\sup_{j=1,\ldots,p}\sup_n\mathbb{E}[\lambda_{j,n}^2]<\infty$;
    \upshape{(b)} $\sup_{j=1,\ldots,p}\sup_n|\lambda_{j,n}|<\infty$ a.s.
\end{lemma}

\begin{lemma}\label{LM_bound_G}
    If Assumptions \ref{A.q}-\ref{A.X-Lipschitz}, \ref{A.Schwartz}, \ref{A.direction},
    \ref{A.beta_uni}, \ref{A.multi_timescale_uni},
    and \ref{A.kernel} hold, then
    \upshape{(a)} $\sup_{n}\mathbb{E}[\|G_{\nu,n}\|^2]<\infty$, \allowbreak
    $\sup_{j=1,\ldots,p}\sup_{n}\mathbb{E}[\|G_{j,n}\|^2]<\infty$;
    \upshape{(b)} $\sup_{n}\|G_{\nu,n}\|<\infty$, $\sup_{j=1,\ldots,p}\sup_{n}\|G_{j,n}\|<\infty$ a.s.
\end{lemma}

\subsection{Asymptotic Behavior of Measure and Gradient Estimates}\label{sec.converge-measure}
Next, we study the asymptotic properties of $\{(\theta_n,\nu_n,\bdnew{\lambda}_n)\}$ through constructing a piecewise continuous-time interpolation of the sequence and its corresponding shifted processes along the timescale characterized by $\{\beta_n\}$. To this end, we let $t_0:=0$ and $t_n:=\sum_{k=0}^{n-1}\beta_{k+1}$ for $n\geq1$. Define $m(t):=\{n:t_n\leq t<t_{n+1}\}$ for $t\geq 0$, and $m(t)=0$ for $t<0$. Let $\nu_\beta^0$ be a continuous-time interpolation of $\{\nu_{n}\}$ constructed as follows: $\nu_\beta^{0}(t):=\nu_{0}$ for $t<0$ and $\nu_\beta^{0}(t):=\nu_{n}$ for $t_n\leq t<t_{n+1}$. We define the family of shifted processes $\{\nu\idxbeta{}\}$ of $\nu_\beta^{0}(\cdot)$ given by $\nu\idxbeta{}(t):=\nu_\beta^{0}(t_n+t)$ for $t\in \Re$. Analogously, we construct the interpolations of $\{(\theta_n,\bdnew{\lambda}_n)\}$ and their shifted processes, denoted by $\{(\theta\idxbeta{},\bdnew{\lambda}\idxbeta{})\}$.
Thus, from \eqref{Eq.nu}
and \eqref{Eq.theta}, the dynamics of $\{(\theta\idxbeta{}(\cdot),\nu\idxbeta{}(\cdot),\lambda\idxbeta{}(\cdot))\}$ can be described by the following equations:
\begin{align*}
    \theta\idxbeta{}(t)=&\theta\idxbeta{}(0)
    +\vartheta\idxbeta{}(t),\\
    \nu\idxbeta{}(t)=&\nu\idxbeta{}(0)
    +\sumbeta\beta_{k+1}\bar{q}(\theta_k, \nu_k)
    +\cm\idxbeta{\nu,}(t),\\
    \lambda\idxbeta{j,}(t)=&\lambda\idxbeta{j,}(0)
    +\cm\idxbeta{j,}(t)
    +\cb\idxbeta{j,}(t)
    +\sumbeta\beta_{k+1}
    f(\nu_k,\theta_k)
    \bar{m}_j(\theta_k,\lambda_{1,k},\ldots,\lambda_{j,k},\nu_k),
\end{align*}
for $j=1,2,\ldots,p$ and $t\ge0$, where we have defined
{\small\begin{align*}
    \vartheta\idxbeta{}(t):=&\sumbeta\alpha_{k+1}(-G_{k}+Z_k),\\
    &\hspace{-16mm}G_{n}:=\nabla_\theta g(\theta_{n},\lambda_{1,n},\ldots,\lambda_{p,n})
    +\sum_{j=1}^p\frac{\partial g(\theta_{n},\lambda_{1,n},\ldots,\lambda_{p,n})}{\partial \lambda_j}G_{j,n},\\
    \cm\idxbeta{\nu,}(t):=&\sumbeta\beta_{k+1}
    [q(\theta_k,X(\theta_k),\nu_k)-\bar{q}(\theta_k, \nu_k)],\\
    \cm\idxbeta{j,}(t):=&\sumbeta\beta_{k+1}
    [K_k m_j(\theta_k,Y(\theta_k),\lambda_{1,k},\ldots,\lambda_{j,k})
    -\ex{K_k m_j(\theta_k,Y(\theta_k),\lambda_{1,k},\ldots,\lambda_{j,k})\vert\F_k}],\\
    \cb\idxbeta{j,}(t):=&\sumbeta\beta_{k+1}
    [\ex{K_k m_j(\theta_k,Y(\theta_k),\lambda_{1,k},\ldots,\lambda_{j,k})\vert\F_k}
    -f(\nu_k,\theta_k)\bar{m}_j(\theta_k,\lambda_{1,k},\ldots,\lambda_{j,k},\nu_k)].
\end{align*}}
and $\alpha_{k+1}Z_k\in-\cc(\theta_{k+1})$ denotes the real vector with the smallest Euclidean norm needed to keep $\theta_{k}-\alpha_{k+1}G_{k}$ from leaving the feasible region $\Theta$.

We show in Lemmas~\ref{LM_Mbeta}-\ref{LM_Thetabeta} below that the processes
$\cm\idxbeta{\nu,}(\cdot),\cm\idxbeta{j,}(\cdot)$,
$\cb\idxbeta{j,}(\cdot)$,
and $\vartheta\idxbeta{}(\cdot)$ are asymptotically negligible as $n\to\infty$. As a result, the asymptotic behavior of the shifted processes $\{(\theta\idxbeta{}(\cdot),\nu\idxbeta{}(\cdot),\lambda\idxbeta{}(\cdot))\}$ can be captured through the following set of coupled ODEs:
\begin{equation}\label{ODE_measure}
\left\{\begin{array}{l}
    \dot{\theta}_\beta(t)=0,\\
    \dot{\nu}_\beta(t)=\bar{q}(\theta_\beta(t),\nu_\beta(t)),\\
    \dot{\lambda}_{j,\beta}(t)=f(\nu_\beta(t),\theta_\beta(t))
    \bar{m}_j(\theta_\beta(t),
    \lambda_{1,\beta}(t),\ldots,\lambda_{j,\beta}(t),
    \nu_\beta(t)),
\end{array}\right.
\end{equation}
for $j=1,2,\ldots,p$ and $t\ge0$, which has a globally asymptotically stable equilibrium taking the form $(\bar{\theta},\nu(\bar{\theta}),\lambda(\bar{\theta}))$ for some $\bar\theta\in\Theta$.
The detailed proofs of lemmas can be found in
supplementary materials.


\begin{lemma}\label{LM_Mbeta}
    Let $T>0$ be given. If Assumptions~\ref{A.q2}, \ref{A.m2},
    \ref{A.beta_uni}-\ref{A.multi_timescale_uni},
    and \ref{A.kernel} hold, then for $j=1,\ldots,p$ and all $t\in[0,T]$,
    $\limn\cm\idxbeta{\nu,}(t)=\limn\cm\idxbeta{j,}(t)=0~~a.s.$
\end{lemma}

\begin{lemma}\label{LM_Bbeta}
    Let $T>0$ be given. If Assumptions~\ref{A.marginal_distribution}, \ref{A.m2}, \ref{A.Schwartz},
    \ref{A.beta_uni}-\ref{A.multi_timescale_uni},
    and \ref{A.kernel} hold, then for $j=1,\ldots,p$ and all $t\in[0,T]$,
    $\limn \cb\idxbeta{j,}(t)=0~~a.s.$
\end{lemma}

\begin{lemma}\label{LM_Thetabeta}
    Let $T>0$ be given. If Assumptions~\ref{A.q}-\ref{A.Schwartz}, \ref{A.direction},
    \ref{A.stepsize_uni},
    and \ref{A.kernel} hold, then for all $t\in[0,T]$,
    $\limn\norm{\vartheta\idxbeta{}(t)}=0~~a.s.$
\end{lemma}

The main convergence result of this section is presented below.
\begin{proposition}\label{PROP_measure}
    If Assumptions~\ref{A.stats}-\ref{A.kernel} hold, then there exists a constant $\bar{\theta}\in\Theta$ such that $ \limn(\theta_n,\nu_n,\bdnew{\lambda}_n)
    =(\bar{\theta},\nu(\bar{\theta}),\lambda(\bar{\theta}))~~a.s.$
\end{proposition}
\proof
    By applying Lemmas~\ref{LM_Mbeta}-\ref{LM_Thetabeta} and Theorem~2.1 in Chapter~5.2 of \cite{kushner2003sa}, we find that for almost every $\omega\in\Omega$, the sequence of shifted processes $\{(\theta\idxbeta{}(\omega,\cdot),\nu\idxbeta{}(\omega,\cdot),\lambda\idxbeta{}(\omega,\cdot))\}$ is equicontinuous in the extended sense on every finite time interval. This, together with the Arzel\`{a}-Ascoli Theorem, implies that the limit $(\theta_\beta(\cdot),\nu_\beta(\cdot),\lambda_\beta(\cdot))$ of any convergent subsequence of the processes
    satisfies the ODEs in \eqref{ODE_measure}.

    For $\theta_\beta(t)\equiv\bar{\theta}\in\Theta$, we have from \eqref{ODE_measure} that
    $$
    \left\{\begin{array}{l}
        \dot{\nu}_\beta(t)=
        \bar{q}\big(\bar{\theta},\nu_\beta(t)\big),\\
        \dot{\lambda}_{j,\beta}(t)=
        f\big(\nu_\beta(t),\bar\theta\big)
        \bar{m}_j\big(\bar\theta,
        \lambda_{1,\beta}(t),\ldots,\lambda_{j,\beta}(t),
        \nu_\beta(t)\big).
    \end{array}\right.
    $$
    By the definition of contextual measures \eqref{def.lam}-\eqref{def.nu}, it is easy to see that the system of ODEs has a unique equilibrium point $(\nu_\beta(t),\lambda_\beta(t))=(\nu(\bar{\theta}),\lambda(\bar{\theta}))$. Let $V(\nu,\lambda)=\ell(\nu-\nu(\bar{\theta}))^2/2+\sum_{j=1}^p\ell_j(\lambda_j-\lambda_j(\bar\theta))^2/2$ be a candidate Lyapunov function for constants $\ell,\ell_1,\ldots,\ell_p>0$. Its time derivative along the ODE trajectory is given by
    \begin{align}
        \dot{V}(\nu_\beta(t),\lambda_\beta(t))
        =&\ell\big(\nu_\beta(t)-\nu(\bar{\theta})\big)
        \bar{q}\big(\bar{\theta},\nu_\beta(t)\big)\label{Eq_dv_sp}\\
        &+f\big(\nu_\beta(t),\bar\theta\big)
        \sum\nolimits_{j=1}^p \ell_j \big(\lambda_{j,\beta}(t)-\lambda_j(\bar{\theta})\big)
        \bar{m}_j\big(\bar\theta,
        \lambda_{1,\beta}(t),\ldots,\lambda_{j,\beta}(t),
        \nu_\beta(t)\big).\label{Eq_dl_sp}
    \end{align}
    From Assumption~\ref{A.q2}, the term~\eqref{Eq_dv_sp} is negative for any $\nu_\beta(t)\neq\nu(\bar{\theta})$.
    On the other hand, given $\nu_\beta(t)=\nu(\bar{\theta})$, we have from Assumption~\ref{A.m2} that the first term in \eqref{Eq_dl_sp}, i.e., $(\lambda_{1,\beta}(t)-\lambda_1(\bar{\theta}))
    \bar{m}_1(\bar\theta,\lambda_{1,\beta}(t),\nu(\bar{\theta}))$, is negative for any $\lambda_{1,\beta}(t)\ne\lambda_1(\bar{\theta})$.
    Next, proceeding iteratively over $j$, we have that for given $(\nu_\beta(t),\lambda_{1,\beta}(t),\ldots,\lambda_{j-1,\beta}(t))=(\nu(\bar{\theta}),\lambda_1(\bar{\theta}),\ldots,\lambda_{j-1}(\bar{\theta}))$, the $j$th term in \eqref{Eq_dl_sp}, i.e., $(\lambda_{j,\beta}(t)-\lambda_j(\bar{\theta}))
    \bar{m}_j(\bar\theta,\lambda_1(\bar{\theta}),\ldots,\lambda_{j-1}(\bar{\theta}),\lambda_{j,\beta}(t),\nu(\bar{\theta}))$, is negative for any $\lambda_{j,\beta}(t)\ne\lambda_j(\bar{\theta})$. This, together with the fact that $\ell,\ell_1,\ldots,\ell_p$ are arbitrary and positive constants, implies that the time derivative $\dot{V}(\nu_\beta(t),\lambda_\beta(t))$ is negative for any $(\nu_\beta(t),\lambda_\beta(t))\neq(\nu(\bar{\theta}),\lambda(\bar{\theta}))$.
    This implies that the equilibrium point $(\bar{\theta},\nu(\bar{\theta}),\lambda(\bar{\theta}))$ is globally asymptotically stable. Finally, because the shifted process
    is constructed by interpolating $\{(\theta_n,\nu_n,\bdnew{\lambda}_n)\}$, we conclude that $\limn(\theta_n,\nu_n,\bdnew{\lambda}_n)=(\bar{\theta},\nu(\bar{\theta}),\lambda(\bar{\theta}))$ a.s.
\endproof

The asymptotic behavior of the sequence of gradient estimates $\{(G_{\nu,n},G_{1,n},\ldots,G_{p,n})\}$ can be analyzed by following a completely analogous argument as above, leading to the following proposition, which indicates the strong consistency of these estimates as the underlying parameter vector $\theta_n$ changes over time. The detailed derivation is provided in
\ref{proof_PROP_gradient}.

\begin{proposition}\label{PROP_gradient}
    If Assumptions~\ref{A.stats}-\ref{A.kernel} hold, then
    \begin{align}\label{converge_Gnu}
    \limn\norm{G_{\nu,n}-\nabla_\theta\nu(\theta)|_{\theta=\theta_n}}
    =\sup_{j=1,\ldots,p}\limn\norm{G_{j,n}-\nabla_\theta\lambda_j(\theta)|_{\theta=\theta_n}}=0
    ~~~~a.s.
    \end{align}
\end{proposition}

\subsection{Main Convergence Theorem}\label{sec.converge-theta}
We have the following result on the convergence of KBSA:
\begin{theorem}\label{TH_converge}
    If Assumptions~\ref{A.stats}-\ref{A.kernel} hold, then the sequence $\{\theta_n\}$ generated by the KBSA algorithm converges a.s. to some limit set of the projected ODE~\eqref{ODE_theta}. In particular, if $\bar{g}$ is strictly convex on $\Theta$, we further have
    \begin{align*}
        \limn\|\theta_n-\theta^*\|=0~~~~a.s.,
    \end{align*}
    where $\theta^*$ denotes the optimal solution to Problem~\eqref{OPT}.
\end{theorem}
\proof
We begin by constructing the shifted processes based on $\{\theta_n\}$ along the timescale characterized by $\{\alpha_n\}$. Let $\tau_0:=0$ and $\tau_n:=\sum_{k=0}^{n-1}\alpha_{k+1}$ for $n\geq1$. Define $\mu(t):=\{n:\tau_n\leq t<\tau_{n+1}\}$ for $t\ge 0$, and $\mu(t)=0$ for $t<0$. Then we construct the piecewise continuous-time interpolation and corresponding shifted processes, denoted as $\{\theta\idxalpha\}$.
Based on \eqref{Eq.theta}, its dynamics can be written as
\begin{align*}
    \theta\idxalpha(t)=\theta\idxalpha(0)
    -\sumalpha\alpha_{k+1}\nabla_\theta\bar{g}(\theta_k)
    +\cb\idxalpha(t)
    +\cz\idxalpha(t),
\end{align*}
for $t\ge0$, where $\cb\idxalpha(t):=\sumalpha\alpha_{k+1}(\nabla_\theta\bar{g}(\theta_k)-G_k)$ and $\cz\idxalpha(t):=\sumalpha\alpha_{k+1}Z_k$.

Note that
\begin{align*}
    \norm{\nabla_\theta\bar{g}(\theta_k)-G_k}
    \le& \norm{\nabla_\theta g(\theta_k,\lambda)|_{\lambda=\lambda(\theta_k)}
    -\nabla_\theta g(\theta_k,\lambda)|_{\lambda=\lambda_k}}
    +\sum\nolimits_{j=1}^p\Big|\frac{\partial g(\theta_k,\lambda_k)}{\partial \lambda_j}\Big|\norm{G_{j,k}-\nabla_\theta\lambda_j(\theta_k)}\\
    &\quad+\sum\nolimits_{j=1}^p\Big|\frac{\partial g(\theta_k,\lambda(\theta_k))}{\partial \lambda_j}
    -\frac{\partial g(\theta_k,\lambda_k)}{\partial \lambda_j}\Big|
    \norm{\nabla_\theta\lambda_j(\theta_k)},
\end{align*}
which tends to $0$ a.s. as $k\to\infty$ as a result of Assumption~\ref{A.g-conti}, Propositions~\ref{PROP_measure} and \ref{PROP_gradient}. Consequently, given $T>0$, we have that for any $t\in[0,T]$, $\|\cb\idxalpha(t)\|
\le\sumalpha\alpha_{k+1}\|\nabla_\theta\bar{g}(\theta_k)-G_k\|\le T\sup_{k\ge n}\|\nabla_\theta\bar{g}(\theta_k)-G_k\|$ a.s.
which tends to $0$ a.s. as $n\to\infty$.

Next, we argue by contradiction and suppose that $\{\cz\idxalpha(\cdot)\}$ is not equicontinuous on $[0, T]$. This contradiction hypothesis implies the existence of sequences $n_k\to\infty$, $\delta_k\to0^+$, and bounded times $\xi_k\in[0,T]$ such that
$\|\cz_\alpha^{n_k}(\xi_k+\delta_k)
-\cz_\alpha^{n_k}(\xi_k)\|\ge\varepsilon$ for some constant $\varepsilon>0$ and all $k\ge0$.
On the other hand, we have from the proof of Lemma~\ref{LM_Thetabeta} that $\sup_n \|Z_n\|<\infty$ a.s. Consequently, for almost every $\omega\in\Omega$, there exists some constant $C_Z>0$ such that
\begin{align*}
    \norm{\cz_\alpha^{n_k}(\omega,\xi_k+\delta_k)
    -\cz_\alpha^{n_k}(\omega,\xi_k)}
    \le\sum\nolimits_{i=\mu(\tau_{n_k}+\xi_k)}^{\mu(\tau_{n_k}+\xi_k+\delta_k)-1}
    \alpha_{i+1}\norm{Z_i(\omega)}
    \le C_Z \sum\nolimits_{i=\mu(\tau_{n_k}+\xi_k)}^{\mu(\tau_{n_k}+\xi_k+\delta_k)-1}\alpha_{i+1},
\end{align*}
which tends to $0$ as $k\to\infty$. This is a contradiction and we conclude that $\{\cz\idxalpha(\cdot)\}$ is a.s. equicontinuous in the extended sense on every finite time interval.

By applying Theorem~2.1 in Chapter~5.2 of \cite{kushner2003sa}, we find that $\{(\theta\idxalpha(\cdot),\cz\idxalpha(\cdot))\}$ is a.s. equicontinuous in the extended sense on every finite time interval, and the limit $(\theta(\cdot),z(\cdot))$ of any convergent subsequence of the processes satisfies the projected ODE~\eqref{ODE_theta}. Hence, on almost every sample path generated by the algorithm, the sequence $\{\theta_n\}$ converges to the limit set of \eqref{ODE_theta}.

In particular, if $\bar{g}$ is strictly convex on $\Theta$, Problem~\eqref{OPT} has a unique solution $\theta^*$ such that $(\theta-\theta^*)^\intercal\nabla_\theta\bar{g}(\theta^*)\ge0$ for every $\theta\in\Theta$. This implies that $\nabla_\theta\bar{g}(\theta^*)\in-\cc(\theta^*)$ and $\theta(t)=\theta^*$ is an equilibrium to ODE~\eqref{ODE_theta}. Let $V(\theta):=\|\theta-\theta^*\|^2/2$ be a candidate Lyapunov function. Because $z(t)\in-\cc(\theta(t))$ and $\bar{g}$ is strictly convex, the time derivative along the trajectory is given by
\begin{align*}
    \dot{V}(\theta(t))
    =(\theta(t)-\theta^*)^\intercal(-\nabla_\theta\bar{g}(\theta)|_{\theta=\theta(t)}+z(t))
    \le-(\theta(t)-\theta^*)^\intercal\nabla_\theta\bar{g}(\theta)|_{\theta=\theta(t)}
    \le \bar{g}(\theta^*)-\bar{g}(\theta(t)),
\end{align*}
which is negative for any $\theta(t)\ne\theta^*$. Thus, we conclude that $\theta^*$ is globally asymptotically stable and thus $\limn\|\theta_n-\theta^*\|=0$ a.s.
\endproof

\section{Rate of Convergence}\label{sec.rate}
In this section, we characterize the convergence rate of KBSA by developing bounds on the MSEs of the sequences of estimates generated. Then we further introduce an acceleration technique based on high-order kernels to improve the convergence rate of the method.

\subsection{Convergence Rate of KBSA}\label{rate_KBSA}
In addition to the conditions used in Section~\ref{sec.converge}, we further make the following assumptions:
\begin{assumption}\label{A.rate}~\
\begin{enumerate}[itemsep=0pt, parsep=0pt,
label=\upshape(\alph*), ref=\theassumption (\alph*)]
    \item\label{A.stepsize2_uni} The sequences
    $\{\alpha_n\}$, $\{\beta_n\}$,
    $\{c_n\}$, and $\{h_n\}$ are taken to be of the form
    $\alpha_n=C_\alpha/(n+n_a)^a$,
    $\beta_n=C_\beta/(n+n_\beta)^b$,
    $c_n=C_c/(n+n_c)^{\mfc}$, and $h_n=C_h/(n+n_h)^{\mfh}$ for some constants
    $C_\alpha, C_\beta, C_c, C_h>0$,
    $n_\alpha, n_\beta, n_c, n_h\ge0$,
    $a, b\in(1/2,1]$, and
    $(\mfc,\mfh)\in\{(\mfc,\mfh):\mfc>0,\mfh>0,\mfc+\mfh<b/2\}$ such that $1/n\le \alpha_n<1$ for all $n\ge N$ and some integer $N>0$.
    \item\label{A.eigen}
    The measures $\nu(\theta), \lambda_1(\theta),\ldots,\lambda_p(\theta)$ are twice continuously differentiable. Let $H(\theta):=\nabla_\theta^2\bar{g}(\theta)$ be the Hessian matrix with the smallest eigenvalue $\rho(\theta)\ge\varepsilon_\rho$ for some constant $\varepsilon_\rho>0$ and all $\theta$ on the line segments between $\theta_n$ and $\theta^*$.
    \item\label{A.int} The ODE~\eqref{ODE_theta} has a unique globally asymptotically stable equilibrium $\theta^*\in\interior{\Theta}$.
    \item\label{A.g-conti2}
    $g$ is twice continuously differentiable with uniformly bounded derivatives.
\end{enumerate}
\end{assumption}

   The condition on the Hessian matrix in Assumption~\ref{A.eigen} is satisfied when the objective function $\bar{g}(\theta)$ is strongly convex on $\Theta$, a standard assumption frequently adopted in the literature analyzing convergence rates of gradient-based algorithms \cite[cf., e.g.,][]{Hu2023, hu2024, wang2025statistical}.
    Assumption~\ref{A.int} requires that the ODE~\eqref{ODE_theta} has a unique equilibrium $\theta^*$ located in the interior of $\Theta$. Under this assumption, we have $\cc(\theta^*)=\{0\}$ and $\theta^*$ satisfies the first-order optimality condition $\nabla_\theta\bar{g}(\theta)|_{\theta=\theta^*}=0$.
    For example, in the context of risk-neutral contextual optimization (Example~\ref{exmp_expect}), a simplified sufficient condition for Assumptions~\ref{A.eigen}-\ref{A.int} is that the cost function $c(\theta,Y)$ is strongly convex in $\theta$, the joint distribution of $(Y,X)$ does not depend on the decision $\theta$, and the optimal solution lies within the interior of the feasible region.

The convergence rates of the MSEs of the contextual measure and their gradient estimates generated by (\ref{Eq.nu}) and (\ref{Eq.Gnu}) are given respectively in Propositions \ref{PROP_rate_measure_uni} and \ref{PROP_rate_gradient_uni} below, based on which an MSE bound on the sequence of parameter vectors $\{\theta_n\}$ is established in Theorem~\ref{TH_rate_uni}. The proofs of these results are given in \ref{proof_rate}.
\begin{proposition}\label{PROP_rate_measure_uni}
    If Assumptions~\ref{A.stats}-\ref{A.rate} hold, then $\{\nu_n\}$ and $\{\bdnew{\lambda}_n\}$ satisfy that
    \begin{enumerate}
    [leftmargin=0em, itemsep=0pt, parsep=0pt, label=]
        \item $\sqrt{\ex{(\nu_n-\nu(\theta_n))^2}}
        =O\big(\frac{\alpha_{n+1}}{\beta_{n+1}}\big)
        +O(\beta_{n+1}^{1/2});$
        \item $\sqrt{\ex{\|\bdnew{\lambda}_n-\lambda(\theta_n)\|^2}}
        =O\big(\frac{\alpha_{n+1}}{\beta_{n+1}}\big)
        +O\big(\frac{\beta_{n+1}^{1/2}}{h_{n+1}^{1/2}}\big)
        +O(h_{n+1}^2)$.
    \end{enumerate}
\end{proposition}

\begin{proposition}\label{PROP_rate_gradient_uni}
    If Assumptions~\ref{A.stats}-\ref{A.rate} hold, then $\{G_{\nu,n}\}$, $\{G_{1,n}\},\ldots,\{G_{p,n}\}$ satisfy that
    \begin{enumerate}
    [leftmargin=0em, itemsep=0pt, parsep=0pt, label=]
        \item $\sqrt{\EE[\|G_{\nu,n}-\nabla_\theta\nu(\theta_n)\|^2]}
        =O\big(\frac{\alpha_{n+1}}{\beta_{n+1}}\big)
        +O\big(\frac{\beta_{n+1}^{1/2}}{c_{n+1}}\big)
        +O(c_{n+1}^2);$
        \item $\sqrt{\EE[\|G_{j,n}-\nabla_\theta\lambda_j(\theta_n)\|^2]}
        =O\big(\frac{\alpha_{n+1}}{\beta_{n+1}}\big)
        +O\big(\frac{\beta_{n+1}^{1/2}}{c_{n+1}h_{n+1}}\big)
        +O(c_{n+1}^2)
        +O(h_{n+1}^2)$,~~
        for $j=1,\ldots,p$.
    \end{enumerate}
\end{proposition}

\begin{theorem}\label{TH_rate_uni}
    If Assumptions~\ref{A.stats}-\ref{A.rate} hold, then $\{\theta_n\}$ generated by the KBSA algorithm satisfies $$\sqrt{\ex{\norm{\theta_n-\theta^*}^2}}=O\bigg(\frac{\alpha_{n+1}}{\beta_{n+1}}\bigg)
        +O\bigg(\frac{\beta_{n+1}^{1/2}}{c_{n+1}h_{n+1}}\bigg)+O(c_{n+1}^2)
        +O(h_{n+1}^2).$$
\end{theorem}

Propositions~\ref{PROP_rate_measure_uni}-\ref{PROP_rate_gradient_uni} and Theorem~\ref{TH_rate_uni} provide practical guidelines for selection of algorithm-specific parameters, including step-, perturbation-sizes, and bandwidths.
Under Assumptions~\ref{A.stepsize_uni} and the specific forms of the algorithm parameters specified in \ref{A.stepsize2_uni}, the optimal decay rate of the performance bound given in Theorem~\ref{TH_rate_uni} can be determined by solving the following optimization problem:
$$\begin{array}{rl}
    \max& \min\{a-b,b/2-\mfc-\mfh,2\mfc,2\mfh\} \\
    \operatorname{s.t.} & 1/2<b<a\le1,\
    \mfc+\mfh<b/2,\ \mfc>0,\ \mfh>0.
\end{array}$$
This allows us to conclude that the optimal upper bound on the convergence rate of $\{\theta_n\}$ is of order $O(n^{-1/5})$, which is attained by setting $a=1$, $b=4/5$, and $\mfc=\mfh=1/10$.

Note that both (\ref{Eq.nu}) and (\ref{Eq.Gnu}) can be used as standalone procedures for contextual measure and sensitivity estimation at a given parameter vector of interest. This is equivalent to setting the step-size $\alpha_n\equiv0$ in KBSA,
in which case, Propositions~\ref{PROP_rate_measure_uni}-\ref{PROP_rate_gradient_uni} imply that the respective bounds on the convergence rates of $\{\nu_n\}$ and $\{\bdnew{\lambda}_n\}$ are of orders $O(n^{-1/2})$ and $O(n^{-2/5})$, which are attained by setting $b=1$ and $\mfh=1/5$. On the other hand, the bound on the convergence rate of $\{G_{\nu,n}\}$ is of order $O(n^{-1/3})$ with $b=1$ and $\mfc=1/6$, whereas the bound on the convergence rate of $\{(G_{1,n},\ldots,G_{p,n})\}$ attains the order $O(n^{-1/4})$ by setting $b=1$ and $\mfc=\mfh=1/8$.

\subsection{Accelerated KBSA}\label{acce_KBSA}
An inspection of Theorem~\ref{TH_rate_uni} suggests that the bias in the kernel estimates typically has an order of $O(h_n^2)$, which may crucially affect the convergence rate of KBSA. To reduce this bias effect, we introduce an acceleration method that uses higher-order kernels (e.g., \citealp{schucany1977improvement,fan1992bias}) to carry out the smoothing.

Suppose the kernel function $\K{\cdot}$ is $\mfr$-times continuously differentiable with uniformly bounded derivatives ($\mfr\ge2$). We define a $\mfr$th order kernel as follows:
\begin{align*}
    \W{u}:=\sum\nolimits_{l=0}^{\mfr-1}
    \begin{pmatrix} \mfr\\l+1 \end{pmatrix}\frac{1}{l!}
    u^l\operatorname{K}^{(l)}\rbrac{u},~~~~~u\in\Re,
\end{align*}
where the superscript $(l)$ denotes the $l$th order derivative.
We show in the result below that when $\W{\cdot}$ is used to carry out
kernel smoothing, the resultant estimation bias can be reduced to the order of $O(h^{\mfr})$.

\begin{proposition}\label{PROP_high_kernel}
    Let $\Phi$ be a function on $\Re$. Suppose that $f(\cdot,\theta)$ and $\Phi(\cdot)$ have $\mfr$-times bounded continuous derivatives. If Assumption~\ref{A.kernel} holds, then
    \begin{align*}
        \frac{1}{h}\ex{\W{\frac{v-X(\theta)}{h}}\Phi(X)}
        -\Phi(v)f(v,\theta)=O(h^{\mfr}).
    \end{align*}
\end{proposition}
\proof
Define $\Psi_\theta(x):=\Phi(x)f(x,\theta)$ for $x\in \cx(\theta)$ and $\Psi_\theta(x):=0$ for $x\notin \cx(\theta)$. By applying a $(\mfr-1)$th order Taylor series expansion of $\Psi_\theta(\cdot)$ around $(v-hu)$, we have
\begin{align*}
    \Psi_\theta(v)
    =&\sum\nolimits_{i=0}^{\mfr-1}\frac{(-h)^i}{i!}
    \Psi_\theta^{(i)}(v-hu)(-u)^i+O(h^{\mfr}).
\end{align*}
By noting that $\int_{-\infty}^{+\infty} \K{u}du=1$ by Assumption~\ref{A.kernel} and applying a change of variable $x=v-hu$, we further have
\begin{align*}
    \Psi_\theta(v)
    =&\int_{-\infty}^{+\infty}
    \sum\nolimits_{i=0}^{\mfr-1}\frac{(-h)^i}{i!}
    \frac{\partial^i}{\partial h^i}
    \Psi_\theta(v-hu)\K{u}du
    +O(h^{\mfr})\\
    =&\int\limits_{\cx(\theta)}\Phi(x)
    \sum\nolimits_{i=0}^{\mfr-1}\frac{(-h)^i}{i!}
    \frac{\partial^i}{\partial h^i}
    \sbrac{\frac{1}{h}\K{\frac{v-x}{h}}}f(x,\theta)dx
    +O(h^{\mfr})\\
    =&\int_{\cx(\theta)}\Phi(x)\frac{1}{h}
    \sum\nolimits_{i=0}^{\mfr-1}
    \sum\nolimits_{l=0}^i
    \begin{pmatrix} i\\l \end{pmatrix}\frac{1}{l!}
    \operatorname{K}_l\rbrac{\frac{v-x}{h}}
    f(x,\theta)dx
    +O(h^{\mfr})\\
    =&\int_{\cx(\theta)}\Phi(x)\frac{1}{h}
    \sum\nolimits_{l=0}^{\mfr-1}
    \begin{pmatrix} \mfr\\l+1 \end{pmatrix}\frac{1}{l!}
    \operatorname{K}_l\rbrac{\frac{v-x}{h}}
    f(x,\theta)dx
    +O(h^{\mfr})\\
    =&\ex{\Phi(X)\frac{1}{h}
    \sum\nolimits_{l=0}^{\mfr-1}
    \begin{pmatrix} \mfr\\l+1 \end{pmatrix}\frac{1}{l!}
    \operatorname{K}_l\rbrac{\frac{v-X}{h}}}
    +O(h^{\mfr}),
\end{align*}
where we have defined $\operatorname{K}_l(u):=u^l\operatorname{K}^{(l)}(u)$. Hence the desired result following directly from the definition of $\W{\cdot}$.
\endproof

As a result of Proposition~\ref{PROP_high_kernel}, the bounds on the convergence rates given in Section~\ref{rate_KBSA} can be further improved by incorporating $\W{\cdot}$ into KBSA.
The results are summarized in Corollary~\ref{COR_accel} below, whose proof follows from the same line of arguments as in Theorem~\ref{TH_rate_uni}.

\begin{corollary}[Accelerated KBSA]\label{COR_accel}
    Suppose that $\W{\cdot}$ is used in lieu of $\K{\cdot}$ in Algorithm~\ref{alg} and the following conditions hold:
    \begin{enumerate}[itemsep=0pt, parsep=0pt,
    label=\upshape(\alph*), ref=\theassumption (\alph*)]
    \item Both $f(\cdot,\theta)$ and $\bar{m}_j(\theta,\lambda_1,\ldots,\lambda_j,\cdot)$ have $\mfr$-times bounded continuous derivatives for every given $(\theta,\lambda)\in\cn\times\Lambda_\varepsilon$ and $j=1,\ldots,p$;
    \item The integral $\int_{-\infty}^{+\infty}
    |u^{2\mfr}\operatorname{K}^{(\mfr)}(u)|du<\infty$.
    \end{enumerate}
    Then under Assumptions~\ref{A.stats}-\ref{A.rate}, the sequences generated by the accelerated KBSA algorithm satisfy
    \begin{enumerate}
    [leftmargin=0em, itemsep=0pt, parsep=0pt, label=]
        \item $\sqrt{\EE[\|{\bdnew{\lambda}}_n
        -\lambda({\theta}_{n})\|^2]}
        =O\big(\frac{\alpha_{n+1}}{\beta_{n+1}}\big)
        +O\big(\frac{\beta_{n+1}^{1/2}}{h_{n+1}^{1/2}}\big)
        +O(h_{n+1}^{\mfr})$;
        \item $\sqrt{\EE[\|{G}_{j,n}
        -\nabla_\theta{\lambda}_j({\theta}_{n})\|^2]}
        =O\big(\frac{\alpha_{n+1}}{\beta_{n+1}}\big)
        +O\big(\frac{\beta_{n+1}^{1/2}}{c_{n+1}h_{n+1}}\big)
        +O(c_{n+1}^2)+O(h_{n+1}^{\mfr})$,~~for $j=1,\ldots,p$;
        \item $\sqrt{\EE[\|{\theta}_n-{\theta}^*\|^2]}
        =O\big(\frac{\alpha_{n+1}}{\beta_{n+1}}\big)
        +O\big(\frac{\beta_{n+1}^{1/2}}{c_{n+1}h_{n+1}}\big)
        +O(c_{n+1}^2)+O(h_{n+1}^{\mfr}).$
    \end{enumerate}
\end{corollary}

In view of Corollary~\ref{COR_accel}, an upper bound on the convergence rate of the mean absolute errors (MAEs) of $\{\theta_n\}$ is given by $O(n^{-\mfr/(4\mfr+2)})$, which is attained by setting $a=1$, $b=(3\mfr+2)/(4\mfr+2)$, $\mfc=\mfr/(8\mfr+4)$ and $\mfh=1/(4\mfr+2)$. In particular, when KBSA is used for estimating contextual measures and sensitivity analysis (i.e., $\alpha_n\equiv 0$), the optimal bound on the MAEs of $\{\bdnew{\lambda}_n\}$ becomes $O(n^{-\mfr/(2\mfr+1)})$, which is achieved by taking $b=1$ and $\mfh=1/(2\mfr+1)$, whereas the optimal bound for $\{G_{j,n}\}$ is $O(n^{-\mfr/(3\mfr+2)})$ under the choice $b=1$, $\mfc=\mfr/[2(3\mfr+2)]$, and $\mfh=1/(3\mfr+2)$. These results reduce to the corresponding ones in Section~\ref{rate_KBSA} when $\mfr=2$.

\section{Multivariate Conditioning Events}\label{sec.extend}
Although the KBSA framework is developed based on contextual measures defined on a single conditioning event, the method can potentially be applied to a wider range of problems under appropriate modifications.
In what follows, we consider a variant of the framework for handling contextual measures involving multiple conditioning events, in which case the covariate becomes a vector $X(\theta)=(X_1(\theta),\ldots,X_{\mfm}(\theta))^\intercal\in\Re^{\mfm}$. We assume that its observed value $\nu(\theta)=(\nu_1(\theta),\ldots,\nu_{\mfm}(\theta))^\intercal$ is given by the unique solution to the system of equations $\ex{q_i(\theta,X_i(\theta),\nu_i(\theta))}=0$ for $i=1,\ldots,{\mfm}$.

The basic idea to extend Algorithm~\ref{alg} to accommodate this multivariate case is to replace the univariate kernel weights $K_{n-1}$ and $K_{n-1}^\pm$ defined in \eqref{term_Kernel_uni} by their multivariate counterparts
\begin{align}\label{kernel_weight}
    K_{n-1}:=\ck\big(H_n^{-1}(\bdnew{\nu}_{n-1}-X(\theta_{n-1}))\big)/\det(H_n),~~~~
    K_{n-1}^\pm:=\ck\big(H_n^{-1}(\bdnew{\nu}_{n-1}^\pm-X(\theta_{n-1}^\pm))\big)/\det(H_n),
\end{align}
where $H_n\in\Re^{\mfm\times \mfm}$ denotes a symmetric and positive definite bandwidth matrix and $\ck:\Re^{\mfm}\to\Re$ is a multivariate kernel function \citep{hardle1997multivariate}. For example, one simple choice of $H_n$ is the diagonal matrix $\operatorname{diag}(h_{1,n},\ldots,h_{\mfm,n})$, and $\ck$ can be taken as the density function of the standard multivariate normal distribution.
In addition, both \eqref{Eq.nu} and \eqref{Eq.Gnu} in Algorithm~\ref{alg} need to be expanded accordingly to also include the recursions for updating the estimates of $\nu(\cdot)$, $\bdnew{\nu}_{n}:=(\nu_{1,n},\ldots,\nu_{\mfm,n})^\intercal$, and their gradients, $(G_{\nu_1,n},\ldots,G_{\nu_{\mfm},n})$. The perturbed values $\bdnew{\nu}_{n-1}^\pm:=(\nu_{1,n-1}^\pm,\ldots,\nu_{\mfm,n-1}^\pm)^\intercal$ in the kernel weights~\eqref{kernel_weight} will then be calculated through $\nu_{i,n-1}^\pm:=\nu_{i,n-1}\pm\bc{n}\Delta_{n}^\intercal G_{\nu_i,n-1}$ for $i=1,\ldots,\mfm$.

The convergence and convergence rate properties of the multivariate KBSA method can be studied by following the same analysis techniques in Sections~\ref{sec.converge}-\ref{sec.rate}. In particular, if we set $H_n=h_n\mathcal{I}_{\mfm}$ where $\mathcal{I}_{\mfm}\in\Re^{\mfm\times \mfm}$ is the identity matrix, the following results can be established:
    \begin{enumerate}
    [leftmargin=0em, itemsep=0pt, parsep=0pt, label=]
        \item $\sqrt{\EE[\|{\bdnew{\lambda}}_n
        -\lambda({\theta}_{n})\|^2]}
        =O\big(\frac{\alpha_{n+1}}{\beta_{n+1}}\big)
        +O\big(\frac{\beta_{n+1}^{1/2}}{h_{n+1}^{\mfm/2}}\big)
        +O(h_{n+1}^2)$;
        \item $\sqrt{\EE[\|{G}_{j,n}
        -\nabla_\theta{\lambda}_j({\theta}_{n})\|^2]}
        =O\big(\frac{\alpha_{n+1}}{\beta_{n+1}}\big)
        +O\big(\frac{\beta_{n+1}^{1/2}}{c_{n+1}h_{n+1}^{\mfm}}\big)
        +O(c_{n+1}^2)+O(h_{n+1}^2)$,~~for $j=1,\ldots,p$;
        \item $\sqrt{\EE[\|{\theta}_n-{\theta}^*\|^2]}
        =O\big(\frac{\alpha_{n+1}}{\beta_{n+1}}\big)
        +O\big(\frac{\beta_{n+1}^{1/2}}{c_{n+1}h_{n+1}^{\mfm}}\big)
        +O(c_{n+1}^2)+O(h_{n+1}^2).$
    \end{enumerate}

As a result, an optimal MAE bound on the convergence rate of $\{\theta_n\}$ is given by $O(n^{-1/(\mfm+4)})$, which is attained by setting $a=1$, $b=(\mfm+3)/(\mfm+4)$, and $\mfc=\mfh=1/[2(\mfm+4)]$. For contextual measure estimation and sensitivity analysis, the optimal bound on the convergence rate of $\{\bdnew{\lambda}_n\}$ is $O(n^{-2/(\mfm+4)})$ under the choice $b=1$ and $\mfh=1/(\mfm+4)$, and the bound on $\{G_{j,n}\}$ becomes $O(n^{-1/(\mfm+3)})$ under the choice $b=1$ and $\mfc=\mfh=1/[2(\mfm+3)]$. These results include those obtained in Section~\ref{rate_KBSA} as a special case when $\mfm=1$.

\section{Numerical Experiments}\label{sec.experiment}
In this section, we carry out simulation experiments to evaluate the performance of KBSA. In Section~\ref{experi_test}, we test the algorithm on a set of artificially created blackbox functions. Then, we apply the method to a nonlinear portfolio sensitivity analysis problem in Section~\ref{experi_portfolio}. Throughout the experiments, Algorithm~\ref{alg} is implemented using the Gaussian kernel function $\K{u}=(2\pi)^{-1/2} e^{-u^2/2}$. Because the performance measures and their gradients of interest vary largely in magnitude, a (possibly) different step-size is used for each of the recursions in \eqref{Eq.nu} and \eqref{Eq.Gnu}. We denote the respective step-sizes used in \eqref{Eq.nu} as $\beta_n,\beta_{1,n},\ldots,\beta_{p,n}$ and those in \eqref{Eq.Gnu} by $\bar\beta_n,\bar\beta_{1,n},\ldots,\bar\beta_{p,n}$. These step-sizes are chosen to have the same (optimal) decay order based on the results in Section~\ref{rate_KBSA}, but possibly with different constant multipliers.

\subsection{Blackbox Test Functions}\label{experi_test}
We consider the following blackbox functions, with dimensions varying from $2$ to $20$ \citep{Hu2023}:
\begin{enumerate}[itemsep=0pt, parsep=0pt, label=Case \arabic*.]\setlength{\itemindent}{2em}
    \item $Y(\theta)=(2.6(\theta_1^2+\theta_2^2)-4.8\theta_1\theta_2)\xi+X(\theta)+X^2(\theta)/2$, for $\Theta=[-2,2]^2$;
    \item $Y(\theta)=(\sum_{i=1}^{10}(\theta_i-i)^2+1)\xi+X(\theta)+X^2(\theta)/2$, for $\theta_i\in[i-1,i+1]$, $i=1,\ldots,10$;
    \item $Y(\theta)=\xi+\sum_{i=1}^{20}(\theta_i-i)\theta_i+X(\theta)+X^2(\theta)/2$, for $\Theta=[-20,20]^{20}$,
\end{enumerate}
where $\xi\sim N(1,1)$, $X(\theta)\sim N(\theta_1,1)$, and $\xi$ is independent of $X(\theta)$. Note that these distributions are primary used for the purpose of comparing with the true optimal solution, and are not known to the algorithm. For each case, we consider the following objective functions:
\begin{enumerate}[itemsep=0pt, parsep=0pt, label=Cost \arabic*.]\setlength{\itemindent}{2em}
    \item $\lambda(\theta)=\ex{Y(\theta)\vert X(\theta)=1}$, with $m(\theta,y,\lambda)=y-\lambda$ and $q(\theta,x,\nu)=1-\nu$;
    \item $\lambda(\theta)=\text{CoVaR}_{0.5,0.6}(\theta)$, with $m(\theta,y,\lambda)=0.6-\I{y\le\lambda}$ and $q(\theta,x,\nu)=0.5-\I{x\le\nu}$.
\end{enumerate}

In the implementation of KBSA, we set the initial values $\theta_0=(1,1)^\intercal$ in Case $1$, $\theta_0=(0.5,1.5,\ldots,9.5)^\intercal$ in Case $2$, and $\theta_0=(0,\ldots,0)^\intercal$ in Case $3$. The initial measure and gradient estimates are set to $\lambda_0=\nu_0=0$, and   $G_{\nu,0}=G_{\lambda,0}=(0,\ldots,0)^\intercal$.
The decreasing rates of the step size, perturbation size, and bandwidth are guided by the results of Theorem~\ref{TH_rate_uni}. Other constants are chosen through trial and error, based on pilot tests conducted prior to larger-scale implementation.
For the experiment with Cost $1$, we set $\alpha_n=1/(n+10^4)$, $\bar\beta_n=\bar\beta_{1,n}=n^{-4/5}$, $h_n=n^{-1/10}$, $c_n=4n^{-1/10}$ for Cases $1$ and $3$, and $c_n=2n^{-1/10}$ for Case $2$. Because $\nu\equiv 1$, we have $\nu_n\equiv1$ and $G_{\nu,n}\equiv(0,\ldots,0)^\intercal$ for all $n$.
For the experiment with Cost $2$, we set  $\alpha_n=1/(n+10^4)$, $\beta_n=\bar\beta_n=2n^{-4/5}$, $\beta_{1,n}=5n^{-4/5}$, $\bar\beta_{1,n}=25n^{-4/5}$, $h_n=n^{-1/10}$, $c_n=0.5n^{-1/10}$ for Case $1$, $c_n=n^{-1/10}$ for Case $2$, and $c_n=8n^{-1/10}$ for Case $3$.
\begin{figure}[h]
    \centering
    {\subfigure[Convergence (Cost 1).
    ]{
        \includegraphics[width=0.42\textwidth]{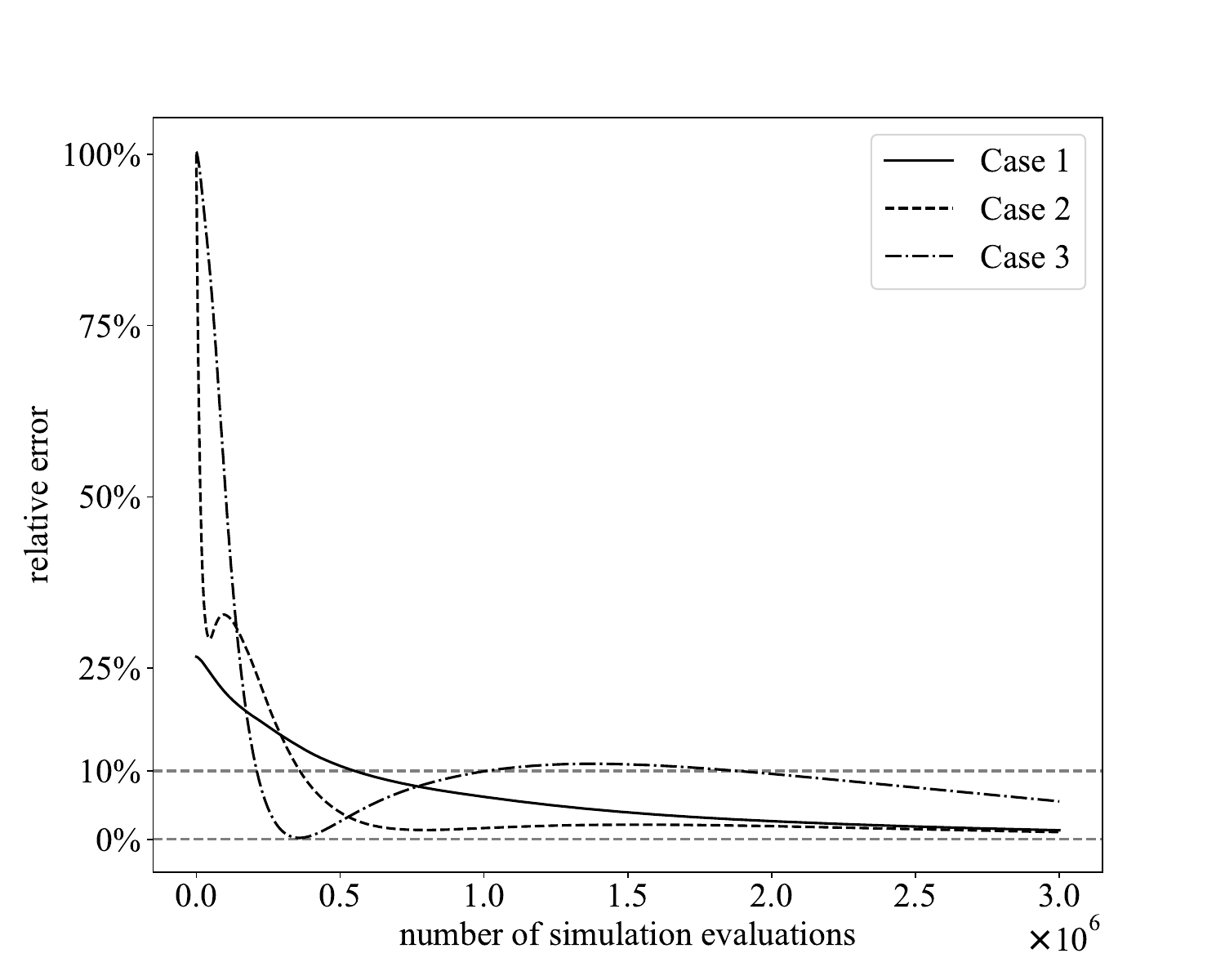}
        \label{fig:converge_mean}
    }
    \subfigure[Rate (Cost 1).]{
        \includegraphics[width=0.42\textwidth]{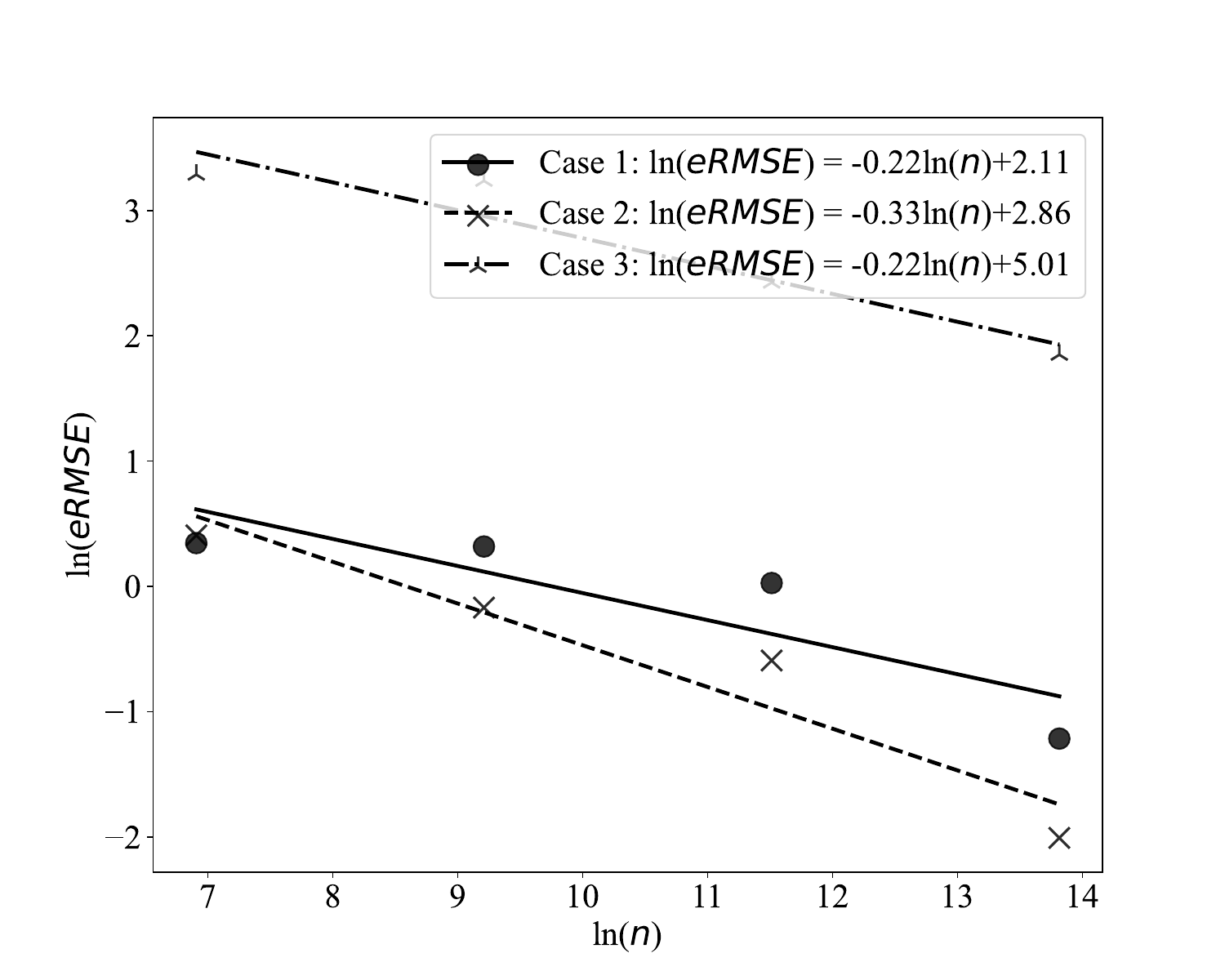}
        \label{fig:rate_test_m}
    }
    \subfigure[Convergence (Cost 2).
    ]{
        \includegraphics[width=0.42\textwidth]{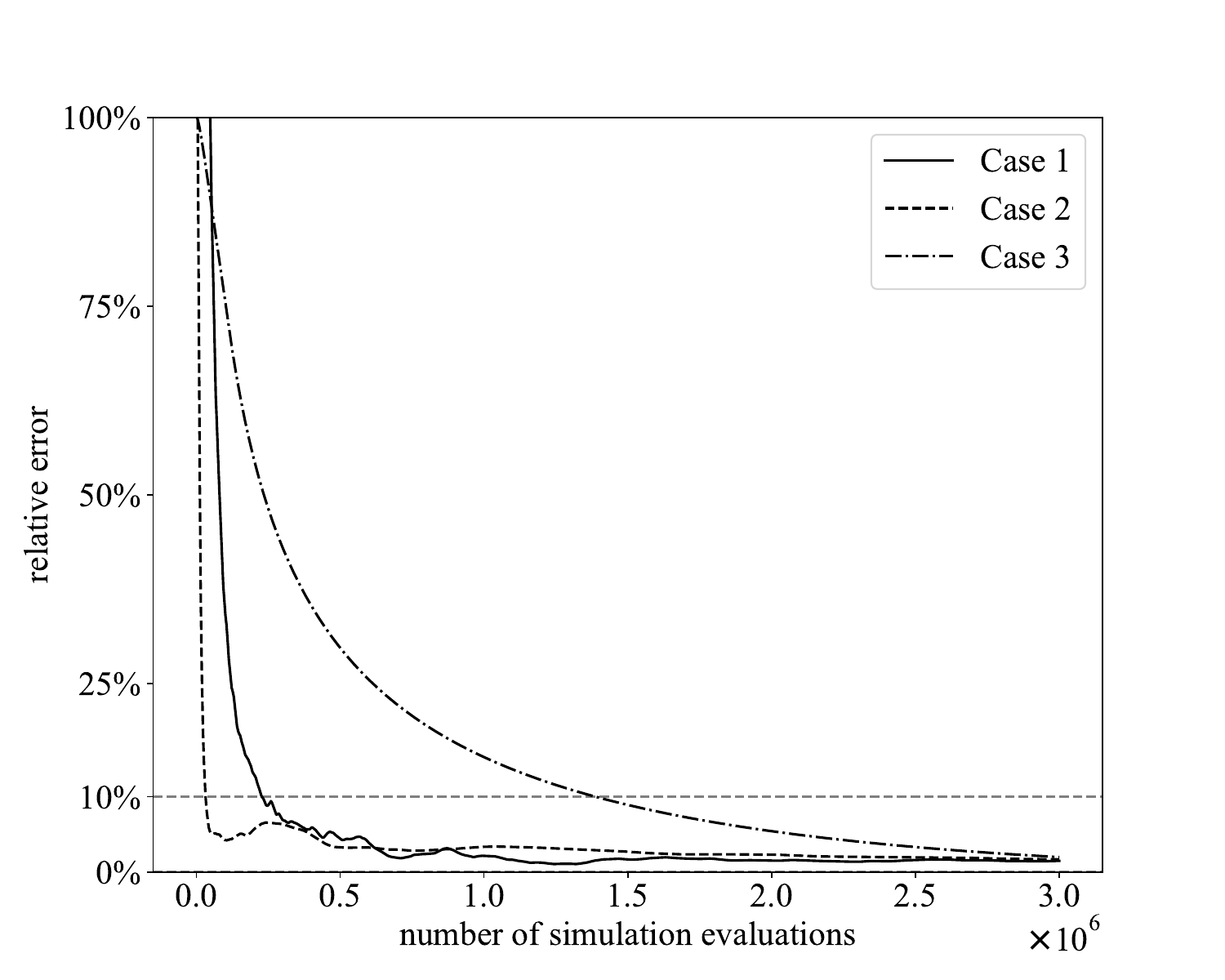}
        \label{fig:converge_risk}
    }
    \subfigure[Rate (Cost 2).]{
        \includegraphics[width=0.42\textwidth]{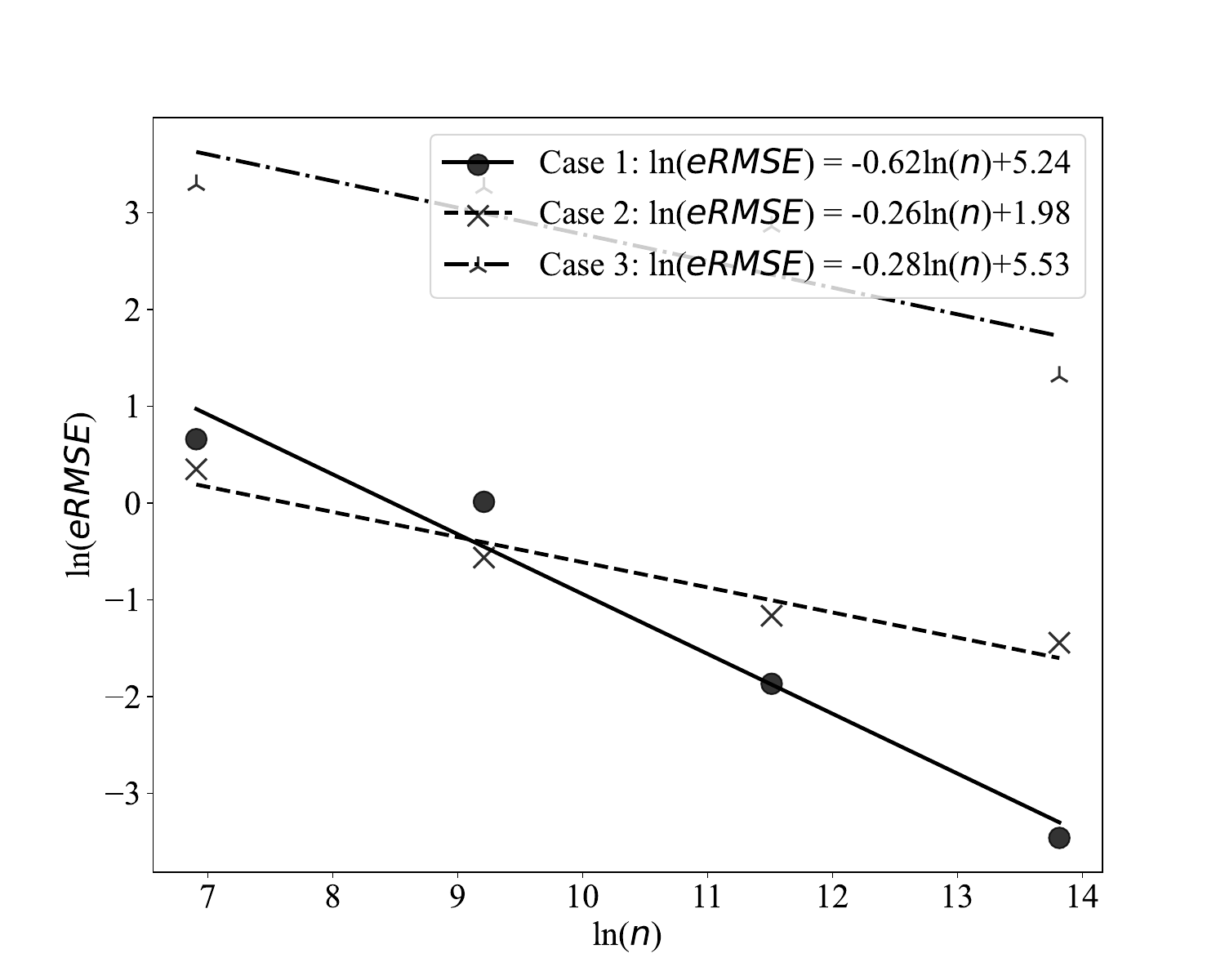}
        \label{fig:rate_test_r}
    }}
    \caption{Convergence behaviors of KBSA (test functions).\label{converge_test}}
\end{figure}

For comparison and algorithm evaluation, we computed the optimal solutions as well as the true values of the measures and gradients in all experiments.
Figure~\ref{fig:converge_mean} and Figure~\ref{fig:converge_risk} illustrate the convergence behavior of KBSA, which plots the relative errors, i.e., $(\text{estimate}-\text{true value})/|\text{true value}|$, of $\{\lambda(\theta_n)\}$ with respect to the optimal value as a function of the number of output evaluations consumed.
Then, experiments are performed with increased numbers of simulation evaluations, $10^3$, $10^4$, $10^5$, and $10^6$. In each case, the algorithm is independently repeated $40$ times. The descriptive statistics of estimators calculated based on the $40$ replications, are reported in Table~\ref{tab:test}, which shows the average relative errors of the objective function values,
i.e., $40^{-1}\sum_{r=1}^{40}[\lambda(\theta_n(r))-\lambda(\theta^*)]/|\lambda(\theta^*)|$,
the empirical root mean squared errors (eRMSEs) of $\{\theta_n\}$,
i.e., $\sqrt{40^{-1}\sum_{r=1}^{40}\|\theta_n(r)-\theta^*\|^2}$, where $\theta_n(r)$ denotes the estimate of $\theta^*$ generated at the $n$th simulation evaluation in the $r$th replication,
and their orders of decrease as the number of iterations $n$ increases.
Figure~\ref{fig:rate_test_m} and Figure~\ref{fig:rate_test_r} also show the log-log plot of the eRMSEs of $\{\theta_n\}$ versus the number of iterations. It can be observed from the figure that the eRMSEs of $\{\theta_n\}$ are of orders $O(n^{-0.22})$, $O(n^{-0.33})$, and $O(n^{-0.22})$ in each of the respective test cases on cost function 1, and of orders $O(n^{-0.62})$, $O(n^{-0.26})$, and $O(n^{-0.28})$ in each of the cases on cost function 2, which conform well to the results obtained in Theorem~\ref{TH_rate_uni}.
In addition, Table~\ref{tab:time_test} reports the average computational times (ACTs) in seconds for the KBSA algorithm across different numbers of iterations.
From Tables~\ref{tab:test} to~\ref{tab:time_test}, we observe that the relative errors of the objective function values fall below $6\%$ once the number of iterations exceeds
$10^6$, while the ACT remains under 10 minutes.
These results demonstrate that KBSA achieves both statistical accuracy and computational efficiency.
\begin{table}[h]
\caption{Relative errors and eRMSEs of estimators on test functions (eRMSEs in parentheses).\label{tab:test}}
\centering\small
{\setlength{\tabcolsep}{2.5mm}
{\begin{tabular}{l|rrrr|l}
\toprule
Cost 1   & $10^3$ & $10^4$        & $10^5$        & $10^6$   &  \\
\hline
Case 1  & {26.63\% (1.41)} & {26.91\% (1.38)}& {14.79\% (1.03)}& {2.09\% (0.30)} & $O(n^{-0.22})$\\
Case 2  & {90.42\% (1.50)} & {28.54\% (0.84)}& {12.26\% (0.55)}& {0.72\% (0.13)} & $O(n^{-0.33})$\\
Case 3  & {100.09\% (26.75)}  & {90.66\% (25.46)}& {17.86\% (11.30)}& {5.65\% (6.35)} & $O(n^{-0.22})$\\
\midrule
Cost 2   & $10^3$ & $10^4$        & $10^5$        & $10^6$   &  \\
\hline
Case 1 & {724.74\% (1.93)}& {199.84\% (1.01)}& {7.60\% (0.15)}& {1.38\% (0.03)} & $O(n^{-0.62})$\\
Case 2 & {114.31\% (1.42)}& {17.18\% (0.57)}& {4.13\% (0.31)}& {1.73\% (0.24)} & $O(n^{-0.26})$\\
Case 3  & {100.03\% (26.76)}& {94.11\% (25.96)}& {42.24\% (17.39)}& {1.91\% (3.69)} & $O(n^{-0.28})$\\
\bottomrule
\end{tabular}}
}
\end{table}
\begin{table}[h]
\caption{Average computational times of KBSA on test functions.\label{tab:time_test}}
\centering\small
{\setlength{\tabcolsep}{2.5mm}
{\begin{tabular}{l|rrrr|l|rrrr}
\toprule
Cost 1   & $10^3$ & $10^4$        & $10^5$        & $10^6$  &
Cost 2   & $10^3$ & $10^4$        & $10^5$        & $10^6$   \\
\hline
Case 1  & {0.19} & {1.90}& {19.18}& {194.09} &
Case 1 & {0.45}& {4.33}& {43.33}& {436.29} \\
Case 2  & {0.24} & {2.43}& {24.39}& {241.61} &
Case 2 & {0.49}& {4.89}& {48.25}& {489.47} \\
Case 3  & {0.28} & {2.90}& {29.24}& {296.67} &
Case 3 & {0.52}& {5.28}& {53.60}& {537.91} \\
\bottomrule
\end{tabular}}
}
\end{table}

\subsection{Nonlinear Portfolios}\label{experi_portfolio}
Consider two portfolios with their losses denoted by $X(\theta)$ and $Y(\theta)$, respectively. Following the assumption used in delta-gamma approximation \citep{glasserman2004monte}, we assume that $X(\theta)\sim N(\theta,\sigma_x^2)$ and
$$
    Y(\theta)=\delta X(\theta)
    +\frac{1}{2}\gamma X^2(\theta)
    +\sigma_y\Big[\rho\Big(\frac{X(\theta)-\theta}{\sigma_x}\Big)
    +\sqrt{1-\rho^2}\xi\Big],
$$
where $\theta=-0.03$, $\sigma_x=0.2$, $\sigma_y=0.3$, $\rho=0.95$, $\delta=0.2$, $\gamma=0.8$, and $\xi$ is a standard normal random variable independent of $X(\theta)$. We are interested in estimating and the sensitivity analysis of the following co-risk measures \citep{Adrian2008, acharya2017measuring}:
\begin{itemize}[itemsep=0pt, parsep=0pt]
    \item $\lambda_1(\theta)=\text{CoVaR}_{0.95,0.95}(\theta)$, with $m_1(\theta,y,\lambda_1)=0.95-\I{y\le\lambda_1}$ and $q(\theta,x,\nu)=0.95-\I{x\le\nu}$;
    \item $\lambda_2(\theta)=\text{CoES}_{0.95,0.95}(\theta)$, with $m_2(\theta,y,\lambda_1,\lambda_2)=(y-\lambda_2)\I{y\ge\lambda_1}$.
\end{itemize}

In the implementation of the proposed algorithm, we set initial estimates $\nu_0=\lambda_{1,0}=\lambda_{2,0}=0$ and $G_{\nu,0}=G_{1,0}=G_{2,0}=1$. Based on the result of Propositions~\ref{PROP_rate_measure_uni}-\ref{PROP_rate_gradient_uni}, we set $\beta_n=
\bar\beta_n=
\beta_{1,n}=\bar\beta_{1,n}=3/(n+1)$,
$\beta_{2,n}=\bar\beta_{2,n}=50/(n+1)$,
$h_n=0.08n^{-1/5}$ for measures estimation, and
$c_n=3n^{-1/8}$, $h_n=0.08n^{-1/8}$ for gradients estimation. To improve the stability of the algorithm, a projection is used in all iterations to project the iterates onto the interval $[-1,1]$.

\begin{figure*}[h]
    \centering
    {\subfigure[Convergence (measure estimates).]{
        \includegraphics[width=0.42\textwidth]{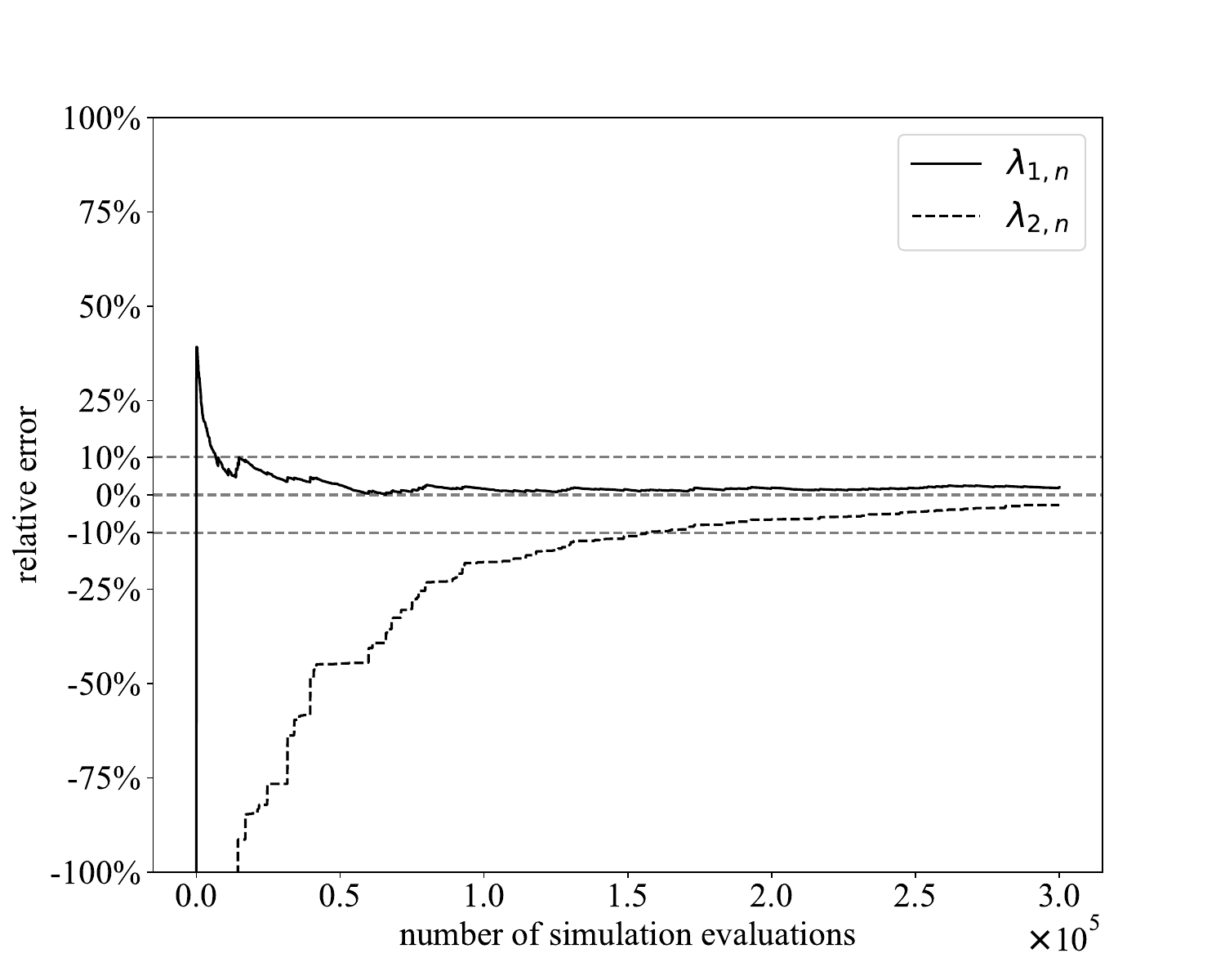}
        \label{fig:converge_port_measure}
    }\subfigure[Rate (measure estimators).]{
        \includegraphics[width=0.42\textwidth]{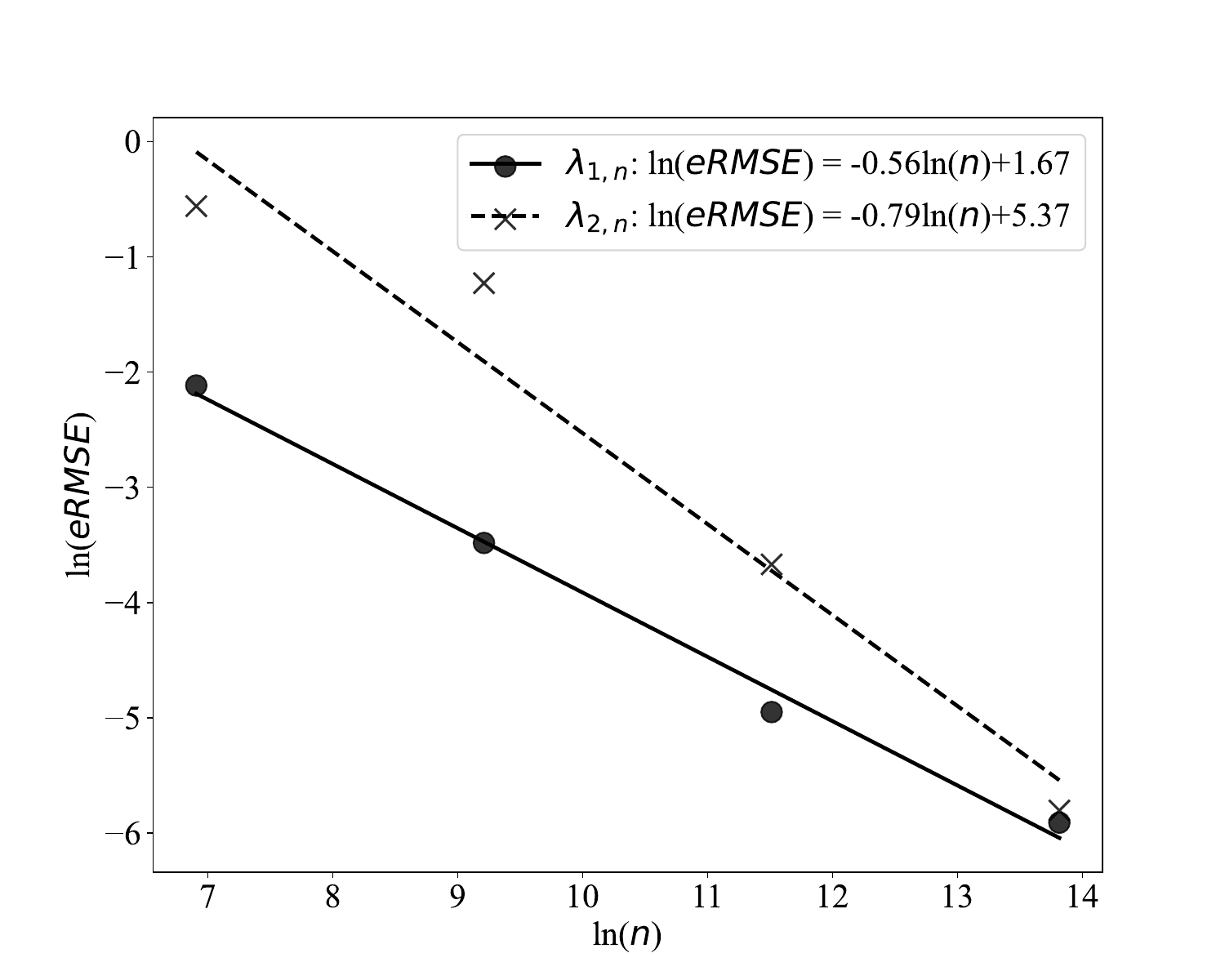}
        \label{fig:rate_port_m}
    }
    \subfigure[Convergence (gradient estimates).]{
        \includegraphics[width=0.42\textwidth]{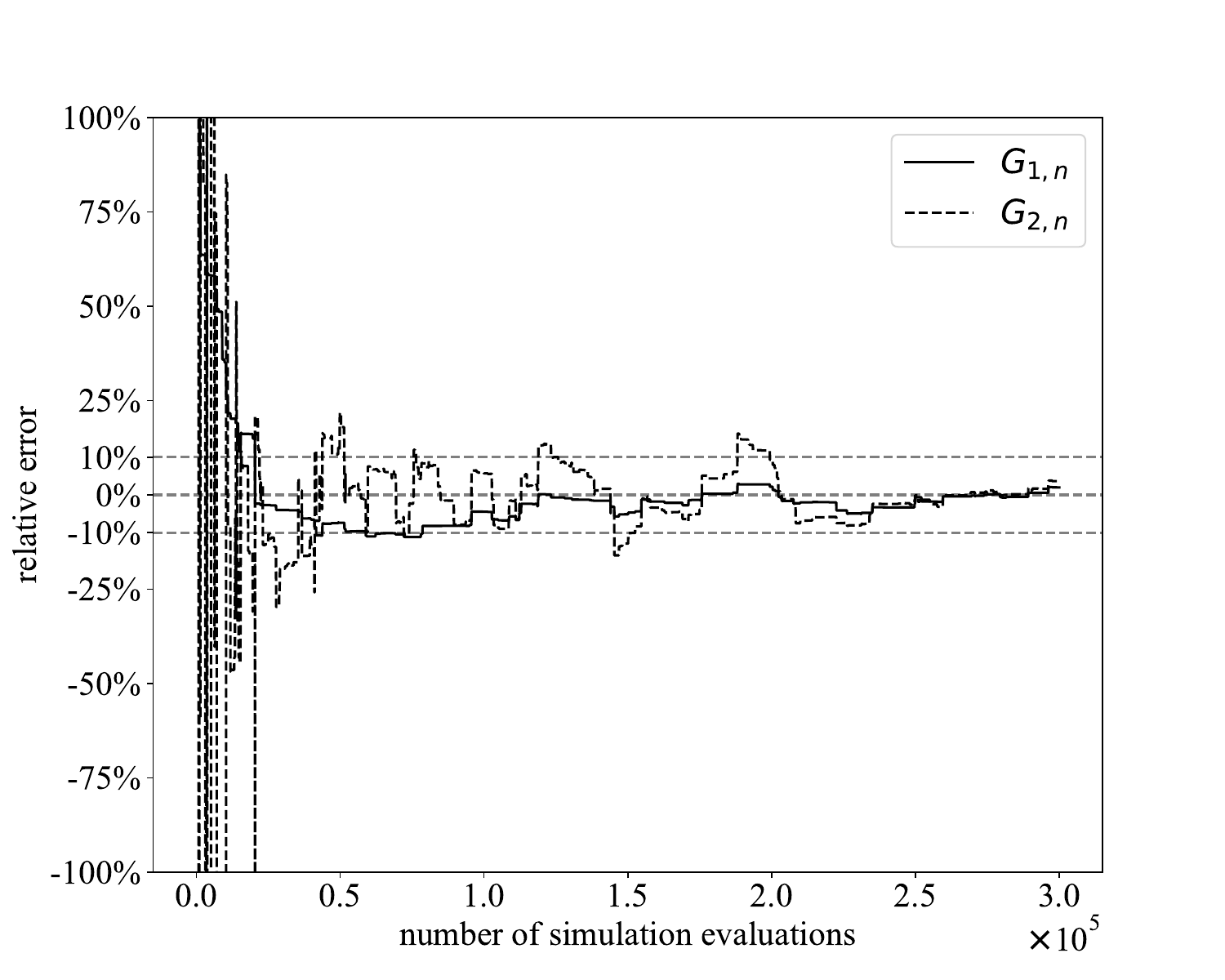}
        \label{fig:converge_port_gradient}
    }\subfigure[Rate (gradient estimators).]{
        \includegraphics[width=0.42\textwidth]{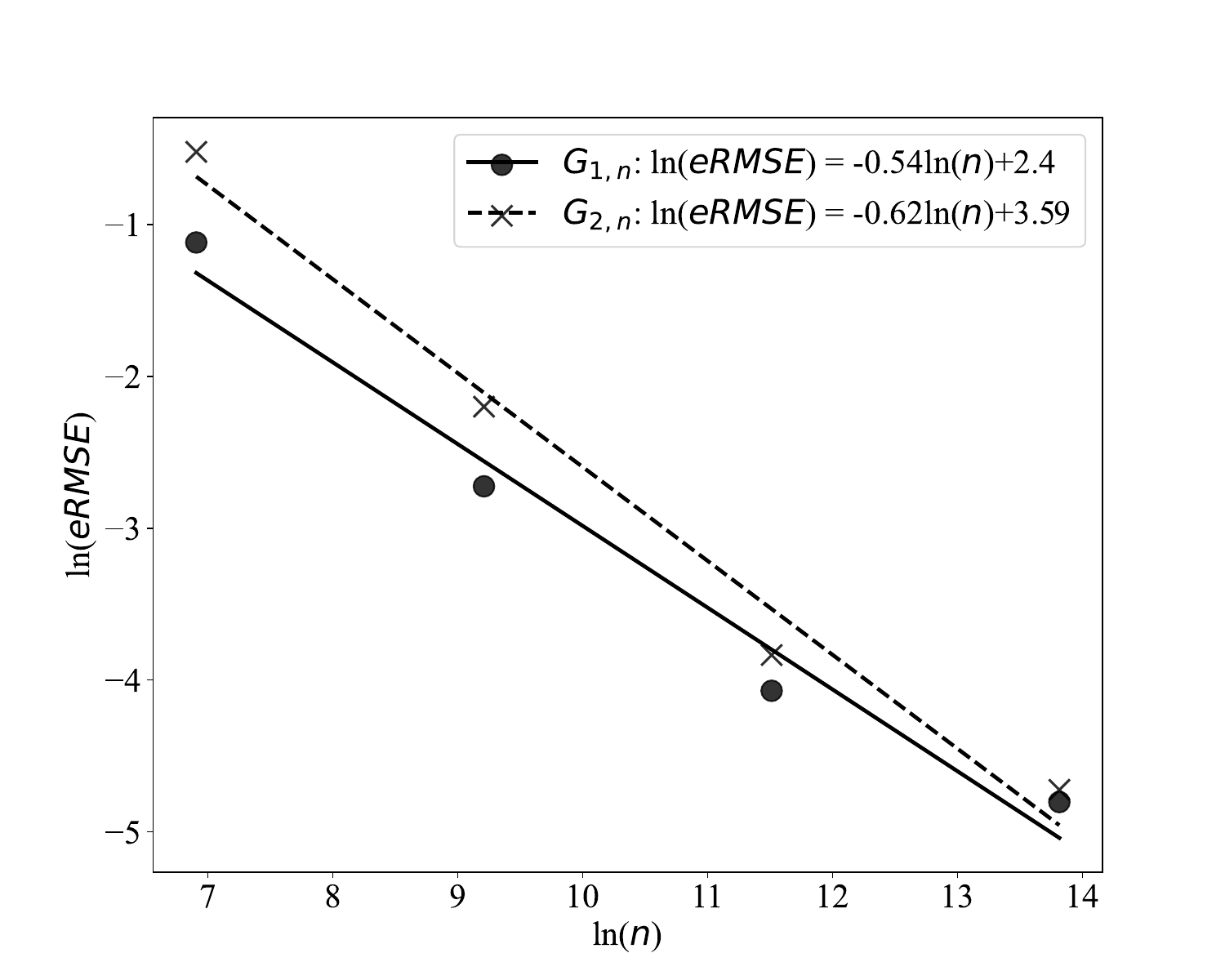}
        \label{fig:rate_port_g}
    }}
    \caption{Convergence behaviors of KBSA (nonlinear portfolios).\label{fig:converge_port}}
\end{figure*}

Figure~\ref{fig:converge_port_measure} and Figure~\ref{fig:converge_port_gradient} show the relative errors of the contextual measure and gradient estimates obtained as functions of the number of iterations. The corresponding eRMSE curves (based on $40$ independent runs) are plotted in Figure~\ref{fig:rate_port_m} and Figure~\ref{fig:rate_port_g}.
The numerical values are reported in Table~\ref{tab:port}. These results clearly indicate the convergence of KBSA and are consistent with the conclusions of Propostions~\ref{PROP_rate_measure_uni}-\ref{PROP_rate_gradient_uni}.

\begin{table}[H]
\caption{Relative errors, eRMSEs, and ACTs of estimators on nonlinear portfolios (eRMSEs in parentheses $\times10^{-2}$).\label{tab:port}}
\centering\small
{\setlength{\tabcolsep}{2.5mm}
{\begin{tabular}{l|rrrr|l}
\toprule
 &  $10^3$   
 & $10^4$        & $10^5$        & $10^6$   &  \\
\hline
$\lambda_{1,n}$    
& {14.04\% (12.05)}& {3.38\% (3.07)}
& {0.45\% (0.71)}& {0.02\% (0.27)} & $O(n^{-0.56})$\\
$\lambda_{2,n}$     
& -33.09\% (57.16)& -14.24\% (29.29)
& -0.31\% (2.55)& 0.10\% (0.30) & $O(n^{-0.79})$\\
\hline
$G_{1,n}$     
& {42.89\% (32.70)}& {7.20\% (6.56)}
& {-0.43\% (1.71)}& {0.06\% (0.82)} & $O(n^{-0.54})$\\
$G_{2,n}$     
& -22.84\% (59.40)& -3.91\% (11.08)
& 0.94\% (2.16)& 0.42\% (0.89) & $O(n^{-0.62})$\\
\midrule
{ACT} & {0.35} 
& {3.65} & {35.57}& {366.29} & 
{$O(n)$} \\
\bottomrule
\end{tabular}}}
\end{table}

\section{Concluding Remarks}\label{sec.conclude}
We have proposed a unifying framework for blackbox optimization of a general class of contextual measures, which covers conditional expectations, conditional quantiles, CoVaR, and CoES. The framework has been extended to analyzing multivariate conditioning events.
Under appropriate conditions, we have shown the strong convergence of KBSA and established the convergence rates of various estimates generated by the framework.
In particular, the convergence rate of KBSA is of order $O(n^{-\mfr/(4\mfr+2)})$ for optimization, whereas when used for contextual measure and gradient estimation, the convergence rates are of orders $O(n^{-\mfr/(2\mfr+1)})$ and $O(n^{-\mfr/(3\mfr+2)})$, respectively, for some integer $\mfr\ge2$.
Simulation studies have also been carried out to illustrate the method, and results conform well with our theoretical findings.

\appendix
\section{Auxiliary Lemma}
We first present Lemma~\ref{LM_kernel_rate}, due to \cite{CHH2023}, that characterizes the approximation error of the kernel method for smooth functions. 
\begin{lemma}\label{LM_kernel_rate}
    Let $\Phi$ be a Schwartz function on $\Re$.  If Assumptions~\ref{A.marginal_distribution}-\ref{A.marginal_distribution2}, \ref{A.Schwartz}, and \ref{A.kernel} hold, then
    \begin{align*}
        \frac{1}{h}\ex{\K{\frac{\nu-X}{h}}\Phi(X)}
        -\Phi(\nu)f(\nu,\theta)
        =\bigo{h^2},~~~~
        \frac{1}{h}\ex{\K{\frac{\nu-X}{h}}
        \rbrac{\Phi(X)-\Phi(\nu)}}
        =&\bigo{h^2},\\
        \ex{\K{\frac{\nu-X}{h}}
        \rbrac{\Phi(X)-\Phi(\nu)}
        }\bigg/
        \ex{\K{\frac{\nu-X}{h}}}
        =&\bigo{h^2}.
    \end{align*}
\end{lemma}

\section{Asymptotic Behavior of the Gradient Estimators}\label{proof_PROP_gradient}
We proceed to analysis of $\{(G_{\nu,n},G_{1,n},\ldots,G_{p,n})\}$. Analogously to Section~\ref{sec.converge-measure}, we construct the continuous-time interpolations and shifted processes of these sequences, denoted as $\{(G\idxbeta{\nu,},G\idxbeta{1,},\ldots,G\idxbeta{p,})\}$.
Besides, we define 
\begin{enumerate}[leftmargin=0em, itemsep=0pt, parsep=0pt, label=]
\item $\db{G_\nu,n}:=\mathbb{E}\Big[\frac{\bar{q}(\theta_{n}^+,\nu_{n}^+)
    -\bar{q}(\theta_{n}^-,\nu_{n}^-)}{2\bc{n+1}\Delta_{n+1}}\Big|\F_n\Big]
    -J(\theta_n,\nu_n,G_{\nu,n}),$
\item $\db{G_j,n}:=\mathbb{E}\Big[K_n^+K_n^-\frac{\bar{m}_{j,n}^+-\bar{m}_{j,n}^-}
    {2\bc{n+1}\Delta_{n+1}}\Big|\F_n\Big]-E_{j,n},$
\item $\tdb{G_j,n}:=E_{j,n}-\bar{E}_j(\theta_n,\nu_n,
    \lambda_{1,n},\ldots,\lambda_{j,n},G_{\nu,n},
    G_{1,n},\ldots,G_{j,n}),$
\item $J(\theta,\nu,G_\nu):=\nabla_\theta\bar{q}(\theta,\nu)
    +G_\nu{\partial \bar{q}(\theta,\nu)}/{\partial\nu},$
\item $\bar{E}_j(\theta,\nu,\lambda_1,\ldots,\lambda_j,
    G_\nu,G_1,\ldots,G_j):=f^2(\nu,\theta)
    \Big(\nabla_\theta\bar{m}_j(\theta,\lambda_1,\ldots,\lambda_j,\nu)$\\
    $\mbox{~~~~}+\sum\nolimits_{i=1}^jG_i{\partial \bar{m}_j(\theta,\lambda_1,\ldots,\lambda_j,\nu)}/{\partial\lambda_i}
    +G_\nu{\partial \bar{m}_j(\theta,\lambda_1,\ldots,\lambda_j,\nu)}/{\partial\nu}\Big),$
\item $Q_{j,h}(\theta,\lambda_1,\ldots,\lambda_j,\nu)
    :={\mathbb{E}[\operatorname{K}((\nu-X(\theta))/h)
    \bar{m}_j(\theta,\lambda_1,\ldots,\lambda_j,X(\theta))]}/{\ex{\operatorname{K}((\nu-X(\theta))/h)}},$
\item $E_{j,n}:=\ex{K_n\vert\F_n}^2\Big(
    \nabla_\theta Q_{j,h_{n+1}}(\theta_n,\lambda_{1,n},\ldots,\lambda_{j,n},\nu_n)$\\
    $\mbox{~~~~}+\sum\nolimits_{i=1}^jG_{i,n}{\partial Q_{j,h_{n+1}}(\theta_n,\lambda_{1,n},\ldots,\lambda_{j,n},\nu_n)}/{\partial\lambda_i}
    +G_{\nu,n}{\partial Q_{j,h_{n+1}}(\theta_n,\lambda_{1,n},\ldots,\lambda_{j,n},\nu_n)}/{\partial\nu}\Big),$
\item $\bar{m}_{j,k}^\pm:=\bar{m}_j(\theta_k^\pm,
    \lambda_{1,k}^\pm,\ldots,\lambda_{j,k}^\pm,X_k^\pm).$
\end{enumerate}

In Lemmas~\ref{LM_Qdiff}-\ref{LM_BK} below, we show that the processes
$\db{G_\nu,n}$, 
$\db{G_j,n}$, 
and $\tdb{G_j,n}$
are asymptotically negligible as $n\to\infty$; the proofs of lemmas are given in supplementary materials. Then we proceed to the proof of Proposition~\ref{PROP_gradient}.

\begin{lemma}\label{LM_Qdiff}
    If Assumptions~\ref{A.marginal_distribution}-\ref{A.marginal_distribution2}, \ref{A.m}, \ref{A.Schwartz}, and \ref{A.kernel} hold, then for $j=1,2,\ldots,p$, all $h\in(0,\delta_h]$, and any fixed $\delta_h>0$,
    \begin{enumerate}[itemsep=0pt, parsep=0pt, 
    label=\upshape(\alph*), ref=\theassumption (\alph*)]
        \item $Q_{j,h}$ is three-times continuously differentiable on $\cn\times\Lambda_\varepsilon\times\cx_\varepsilon$ with uniformly bounded derivatives;
        \item for $i=1,\ldots,j$ and $(\theta,\lambda,\nu)\in
        \cn\times\Lambda_\varepsilon\times\cx_\varepsilon$,
       {\begin{align*}
            \hspace{-4mm}\norm{\nabla_\theta Q_{j,h}(\theta,\lambda_1,\ldots,\lambda_j,\nu)
            -\nabla_\theta \bar{m}_j(\theta,\lambda_1,\ldots,\lambda_j,\nu)}
            =&O(h^2),\\
            \hspace{-4mm}\abs{{\partial Q_{j,h}(\theta,\lambda_1,\ldots,\lambda_j,\nu)}/{\partial\lambda_i} 
            -{\partial \bar{m}_j(\theta,\lambda_1,\ldots,\lambda_j,\nu)}/{\partial\lambda_i}}
            =&O(h^2),\\
            \hspace{-4mm}\abs{{\partial Q_{j,h}(\theta,\lambda_1,\ldots,\lambda_j,\nu)}/{\partial\nu} 
            -{\partial \bar{m}_j(\theta,\lambda_1,\ldots,\lambda_j,\nu)}/{\partial\nu}}
            =&O(h^2).
        \end{align*}}
    \end{enumerate}
\end{lemma}

\begin{lemma}\label{LM_Bgamma}
    If Assumptions~\ref{A.marginal_distribution}-\ref{A.X-Lipschitz}, \ref{A.direction}, \ref{A.beta_uni}, \ref{A.multi_timescale_uni}, and \ref{A.kernel} hold, then 
    $\limn\|\db{G_\nu,n}\|=\sup_j\limn\|\db{G_j,n}\|=0~~a.s.$
\end{lemma}

\begin{lemma}\label{LM_BK}
    If Assumptions~\ref{A.marginal_distribution}-\ref{A.X-Lipschitz}, \ref{A.Schwartz}, \ref{A.direction}, \ref{A.beta_uni}, \ref{A.multi_timescale_uni}, and \ref{A.kernel} hold, then for $j=1,\ldots,p$,
    $\limn\|\tdb{G_j,n}\|=0~~a.s.$
\end{lemma}

\proof[Proof of Proposition~\ref{PROP_gradient}]
    For almost every $\omega\in\Omega$, the sequence of shifted processes 
    $\{(\theta\idxbeta{}(\omega,\cdot),
    \nu\idxbeta{}(\omega,\cdot),
    \lambda\idxbeta{}(\omega,\cdot),\allowbreak
    G\idxbeta{\nu,}(\omega,\cdot),
    G\idxbeta{1,}(\omega,\cdot),\ldots,
    G\idxbeta{p,}(\omega,\cdot))\}$ is equicontinuous on every finite time interval. 
    And the limit of any convergent subsequence of the processes satisfies ODEs~\eqref{ODE_measure} and the following ODEs:
    $$\left\{
    \begin{array}{l}
        \dot{G}_{\nu,\beta}(t)=
        J\big(\theta_\beta(t),\nu_\beta(t),
        G_{\nu,\beta}(t)\big),\\
        \dot{G}_{j,\beta}(t)=
        \bar{E}_j\big(\theta_\beta(t),\nu_\beta(t),
        \lambda_{1,\beta}(t),\ldots,\lambda_{j,\beta}(t),
        G_{\nu,\beta}(t),
        G_{1,\beta}(t),\ldots,G_{j,\beta}(t)\big),
    \end{array}\right.
    $$
    for $j=1,\ldots,p,$ with a globally asymptotically stable equilibrium $
    (\bar\theta,
    \nu(\bar\theta),
    \lambda(\bar\theta),
    \nabla_\theta\nu(\bar\theta),
    \nabla_\theta\lambda_1(\bar\theta),
    \ldots,
    \nabla_\theta\lambda_p(\bar\theta))$. 

    Thus, we find that $\limn\|G_{\nu,n}-\nabla_\theta\nu(\bar\theta)\|=\limn\|G_{j,n}-\nabla_\theta\lambda_j(\bar\theta)\|=0$ a.s., for $j=1,\ldots,p$.
    Finally, note that $\|G_{\nu,n}-\nabla_\theta\nu(\theta)|_{\theta=\theta_n}\|
        \le\|G_{\nu,n}-\nabla_\theta\nu(\bar\theta)\|
        +\|\nabla_\theta\nu(\theta)|_{\theta=\theta_n}-\nabla_\theta\nu(\bar\theta)\|$ and 
        $\|G_{j,n}-\nabla_\theta\lambda_j(\theta)|_{\theta=\theta_n}\|
        \le\|G_{j,n}-\nabla_\theta\lambda_j(\bar\theta)\|
        +\|\nabla_\theta\lambda_j(\theta)|_{\theta=\theta_n}-\nabla_\theta\lambda_j(\bar\theta)\|$ a.s.
    Because $\limn\theta_n=\bar\theta$ a.s. by Proposition~\ref{PROP_measure} and $\nabla_\theta\nu(\cdot)$, $\nabla_\theta\lambda_j(\cdot)$ are continuous by Assumptions~\ref{A.q} and \ref{A.m}, we obtain \eqref{converge_Gnu} as required. 
\endproof

\section{Convergence Rate Analysis}\label{proof_rate}
We first present Lemma~\ref{LM_contraction_mapping}, due to \cite{CHH2023}, that characterizes the order of the iterates produced from a sequence of contraction mappings. Then we proceed to the proof of Propositions~\ref{PROP_rate_measure_uni}-\ref{PROP_rate_gradient_uni} and Theorem~\ref{TH_rate_uni}.
\begin{lemma}\label{LM_contraction_mapping}
    Let $\{u_n\}$, $\{v_n\}$, and $\{w_n\}$ be sequences of positive real numbers. Define the mappings
    $
    \Gamma_n(x):=\sqrt{((1-u_n)x+w_n)^2+v_n},~~x\geq 0,
    $
    for $n\ge1$ and let $\{x_n\}$ be the sequence generated by $x_{n}=\Gamma_n(x_{n-1})$.  Suppose that $u_n=\bar{u}/(n+n_0)^u$ for constants $\bar{u}>0,n_0\ge0,u\in(0,1]$ such that $1/n\le u_n<1$ for all $n\ge N$ and some integer $N>0$, and there exists a constant $s\in(u,u+1)$ such that $|
        {w_{n}}/{u_{n}}
        +\sqrt{{v_{n}}/{u_{n}}}
        -{w_{n+1}}/{u_{n+1}}
        -\sqrt{{v_{n+1}}/{u_{n+1}}}
        |=\bigo{n^{-s}}$ and ${w_{n+1}}/{u_{n+1}}
        +\sqrt{{v_{n+1}}/{u_{n+1}}}
        =\sfomega(n^{u-s})$.
    Then
    $x_n=\bigo{{w_{n+1}}/{u_{n+1}}}
    +\bigo{\sqrt{{v_{n+1}}/{u_{n+1}}}}.$
\end{lemma}
\proof[Proof of Proposition~\ref{PROP_rate_measure_uni}]
    Define $\cv_n:=\nu_n-\nu(\theta_n)$. 
    Because $\nu$ is continuously differentiable and $\Theta$ is compact, we note that $\nu$ is Lipschitz continuous on $\Theta$. This, together with Equation~\eqref{Eq.theta} and the argument used in the proof of Lemma~\ref{LM_Thetabeta}, implies $|\nu(\theta_{n-1})-\nu(\theta_n)|\le L_\nu\|\theta_{n-1}-\theta_n\|\le \alpha_n L_\nu(\norm{Z_{n-1}}+\norm{G_{n-1}})\le 2\alpha_n L_\nu\norm{G_{n-1}}$ a.s.
    for some constant $L_\nu>0$. This, together with Equation~\eqref{Eq.nu}, the mean value theorem, and Assumption~\ref{A.q2}, implies that $\ex{\cv_n^2\vert\F_{n-1}}\le(1-2\beta_n\varepsilon_q)\cv_{n-1}^2
        +4\alpha_nL_\nu\abs{\cv_{n-1}}\norm{G_{n-1}}
        +\beta_n^2C_q
        +4\alpha_n^2L_\nu^2\norm{G_{n-1}}^2
        +4\beta_n\alpha_nL_\nu\sqrt{C_q}\norm{G_{n-1}}~~a.s.$
    
    Let $N:=\inf\{n:\beta_n\varepsilon_q\le1/2\}$. By taking expectations at both sides and applying \holder,  Lemma~\ref{LM_bound_G}, and Assumption~\ref{A.multi_timescale_uni}, we then find that $\mathbb{E}[\cv_n^2]\le((1-\beta_n\varepsilon_q)
        \sqrt{\mathbb{E}[\cv_{n-1}^2]}+\alpha_n\bar{C}_\nu)^2
        +\beta_n^2C_1$
    for all $n\ge N$ and some constants $C_G,C_1>0,\bar{C}_\nu:=2L_\nu \sqrt{C_G}/\sqrt{1-2\beta_N\varepsilon_q}$.

    Consider the sequence $\{x_n\}$ generated by $x_n=[((1-\beta_n\varepsilon_q)
    x_{n-1}+\alpha_n\bar{C}_\nu)^2+\beta_n^2C_1]^{1/2}$
    for $n\ge N$, and $x_{N-1}:=\sqrt{\mathbb{E}[\cv_{N-1}^2]}$. By applying an inductive argument, we have $\sqrt{\mathbb{E}[\cv_n^2]}\le x_n$ for all $n\ge N$. On the other hand, we have from Assumptions~\ref{A.multi_timescale_uni} and \ref{A.stepsize_uni} that
    $|{\alpha_n}/{(\beta_n\varepsilon_q)}
    +\sqrt{{\beta_nC_1}/{\varepsilon_q}}
    -{\alpha_{n+1}}/{(\beta_{n+1}\varepsilon_q)}
    -\sqrt{{\beta_{n+1}C_1}/{\varepsilon_q}}|=O(n^{-s})$ and ${\alpha_{n+1}}/{(\beta_{n+1}\varepsilon_q)}
    +\sqrt{{\beta_{n+1}C_1}/{\varepsilon_q}}=\sfomega(n^{u-s}),$
    for $s=\min\{a-b,b/2\}+1$ and $u=s-\max\{a-b,b/2\}$. Thus, the proof of part (a) is completed by applying Lemma~\ref{LM_contraction_mapping}.

    On the other hand, define $\cl_{j,n}:=\lambda_{j,n}-\lambda_j(\theta_n)$. Analogously, we have from Equation~\eqref{Eq.nu},
    Lemma~\ref{LM_kernel_rate}, 
    the argument used in the proof of Lemma~\ref{LM_bound_lam},
    and the fact that $\abs{\lambda_j(\theta_{n-1})-\lambda_j(\theta_n)}\le 2\alpha_n L_j\norm{G_{n-1}}$ for some constant $L_j>0$, that $\mathbb{E}[\cl_{1,n}^2\vert\F_{n-1}]\le (1-2\beta_{n}\varepsilon_m\varepsilon_f)\cl_{1,n-1}^2
        +\beta_{n}\bar{C}_1|\cl_{1,n-1}|\abs{\cv_{n-1}}
        +\beta_{n}h_n^2\bar{C}_2|\cl_{1,n-1}|
        +4\alpha_nL_1|\cl_{1,n-1}|\norm{G_{n-1}}
        +4\alpha_n^2 L_1^2\norm{G_{n-1}}^2
        +\beta_{n}^2C_KC_f\allowbreak C_m/{h_n}
        +4\alpha_n\beta_{n}L_1C_f\sqrt{C_m}
        \norm{G_{n-1}}$ a.s.,
    for some constants $\bar{C}_1,\bar{C}_2>0$. Let $\bar{N}:=\inf\{n:\beta_{n}\varepsilon_m\varepsilon_f\le1/2\}$. By taking expectations and applying the same argument used in the proof of part (a), we then find that $\mathbb{E}[\cl_{1,n}^2]\le((1-\beta_{n}\varepsilon_m\varepsilon_f)
        \sqrt{\mathbb{E}[\cl_{1,n-1}^2]}+\ca_{1,n})^2
        +{(\beta_{n}^2\tilde{C}_3)}/{h_n},$
    for all $n\ge\bar{N}$ and some constants $\tilde{C}_1,\tilde{C}_2,\tilde{C}_3>0$, where $\ca_{1,n}:=(1-2\beta_{\bar{N}}\varepsilon_m\varepsilon_f)^{-1/2}[(\beta_n/2)(\alpha_n\tilde{C}_1/\beta_n+\beta_n^{1/2}\tilde{C}_2+h_n^2\bar{C}_2)
    +2\alpha_nL_1\sqrt{C_G}]$.
    Consequently, by applying Lemma~\ref{LM_contraction_mapping} and essentially the same argument used in the proof of part (a), we find that $\sqrt{\mathbb{E}[\cl_{1,n}^2]}
    =O({\alpha_{n+1}}/{\beta_{n+1}})
    +O(\beta_{n+1}^{1/2}/h_{n+1}^{1/2})+O(h_{n+1}^2)$. Analogously, it is straightforward to show that $\sqrt{\mathbb{E}[\cl_{j,n}^2]}
    =O({\alpha_{n+1}}/{\beta_{n+1}})
    +O(\beta_{n+1}^{1/2}/h_{n+1}^{1/2})+O(h_{n+1}^2)$ by induction and essentially the same argument.
\endproof

The proof of Proposition~\ref{PROP_rate_gradient_uni} follows an analogous argument as above and thus is given in supplementary materials.

\proof[Proof of Theorem~\ref{TH_rate_uni}]
    Define $\zeta_n:=\theta_n-\theta^*$. By applying Equation~\eqref{Eq.theta}, \holder,  and the argument used in the proof of Lemma~\ref{LM_Thetabeta}, we find that $\norm{\zeta_n}^2\le\norm{\zeta_{n-1}}^2-2\alpha_n\zeta_{n-1}^\intercal
        [\nabla_\theta\bar{g}(\theta_{n-1})-\nabla_\theta\bar{g}(\theta^*)]
        +2\alpha_n\norm{\zeta_{n-1}}\norm{\cg_{n-1}}
        +2\alpha_n\norm{\zeta_{n-1}}\norm{Z_{n-1}}
        +4\alpha_n^2\norm{G_{n-1}}^2~~a.s.,$
    where $\cg_{n-1}:=G_{n-1}-\nabla_\theta\bar{g}(\theta_{n-1})$, $\bar\theta_{n-1}$ is on the line segment between $\theta_{n-1}$ and $\theta^*$. Then by applying the mean value theorem, Rayleigh-Ritz inequality, and Assumption~\ref{A.eigen}, we find that
    \begin{align}
        \norm{\zeta_n}^2
        \le&\norm{\zeta_{n-1}}^2-2\alpha_n\zeta_{n-1}^\intercal H(\bar\theta_{n-1})\zeta_{n-1}
        +2\alpha_n\norm{\zeta_{n-1}}\norm{\cg_{n-1}}
        +2\alpha_n\norm{\zeta_{n-1}}\norm{Z_{n-1}}
        +4\alpha_n^2\norm{G_{n-1}}^2\nonumber\\
        \le&(1-2\alpha_n\varepsilon_\rho)\norm{\zeta_{n-1}}^2
        +2\alpha_n\norm{\zeta_{n-1}}\norm{\cg_{n-1}}
        +2\alpha_n\norm{\zeta_{n-1}}\norm{Z_{n-1}}
        +4\alpha_n^2\norm{G_{n-1}}^2~~~~a.s.,\label{tmp_zeta2}
    \end{align}
    for all $n\ge N$ and $N:=\inf\{n:\alpha_n\varepsilon_\rho\le1/2\}$, where we have from Assumption~\ref{A.g-conti2} that $\norm{\cg_{n-1}}\le\sum_{j=1}^pC_{\lambda,j}\abs{\cl_{j,n-1}}
    +\sum_{j=1}^pC_{G,j}\norm{\cg_{j,n-1}}$ a.s.
    for some constants $C_{\lambda,1},\ldots,C_{\lambda,p},C_{G,1},\ldots,C_{G,p}>0$, $\cl_{j,n-1}:=\lambda_{j,n-1}-\lambda_j(\theta_{n-1})$, and $\cg_{j,n-1}:=G_{j,n-1}-\nabla_\theta\lambda_j(\theta_{n-1})$. 

    Because $\theta^*\in\interior{\Theta}$, there is a constant $\delta>0$ such that $\scru_\delta(\theta^*):=\{\vartheta\in\Re^d:\|\vartheta-\theta^*\|\le2\delta\}\subseteq\Theta$. By applying Markov's inequality, we find that
    \begin{align}
        \ex{\norm{Z_{n-1}}^2}
        =&\ex{\norm{Z_{n-1}}^2\vert\theta_n\in\scru_\delta(\theta^*)}
        \pr{\theta_n\in\scru_\delta(\theta^*)}
        +\ex{\norm{Z_{n-1}}^2\vert\theta_n\notin\scru_\delta(\theta^*)}
        \pr{\theta_n\notin\scru_\delta(\theta^*)}\nonumber\\
        \le&\delta^{-2}({4\alpha_n^2\mathbb{E}[\norm{G_{n-1}}^2]^2}
        +{\mathbb{E}[\norm{G_{n-1}}^2]\mathbb{E}[\norm{\zeta_{n-1}}^2}]),\label{tmp_Z2}
    \end{align}
    where we have $\limn\mathbb{E}[\norm{G_n}^2]=0$ because $\norm{G_n}\le\norm{\cg_n}+\norm{\nabla_\theta\bar{g}(\theta_n)}$, $\limn\mathbb{E}[\norm{\cg_n}^2]=0$, and $\limn\mathbb{E}[\norm{\nabla_\theta\bar{g}(\theta_n)}^2]=0$.

    As a result, for given $\varepsilon_G\in(0,\delta^2\varepsilon_\rho^2)$, there exists an integer $N_G$ such that $\mathbb{E}[\norm{G_n}^2]\le\varepsilon_G$ for all $n\ge N_G$. Let $\bar{N}:=\inf\{n\ge N_G:\alpha_n(\varepsilon_\rho-\sqrt{\varepsilon_G}/\delta)\le1/2\}$. 
    By taking expectations at both sides of Equation~\eqref{tmp_zeta2} and applying \holder, Equation~\eqref{tmp_Z2}, Assumption~\ref{A.multi_timescale_uni}, and Propositions~\ref{PROP_rate_measure_uni}-\ref{PROP_rate_gradient_uni}, we find that
    \begin{align*}
        \ex{\norm{\zeta_n}^2}
        \le&(1-2\alpha_n\varepsilon_\rho)\ex{\norm{\zeta_{n-1}}^2}
        +4\alpha_n^2\varepsilon_G
        +2\alpha_n\sqrt{\ex{\norm{\zeta_{n-1}}^2}}
        \sum\nolimits_{j=1}^pC_{\lambda,j}\sqrt{\ex{\cl_{j,n-1}^2}}\\
        &\quad+2\alpha_n\sqrt{\ex{\norm{\zeta_{n-1}}^2}}
        \sum\nolimits_{j=1}^pC_{G,j}\sqrt{\ex{\norm{\cg_{j,n-1}}^2}}
        +2(\alpha_n/\delta)\sqrt{\ex{\norm{\zeta_{n-1}}^2}}
        ({2\alpha_n\varepsilon_G}
        +{\sqrt{\varepsilon_G\ex{\norm{\zeta_{n-1}}^2}}})\\
       \le&\Big((1-\alpha_n\varepsilon_\zeta)\sqrt{\ex{\norm{\zeta_{n-1}}^2}}
        +\ca_{\zeta,n}\Big)^2+4\alpha_n^2\varepsilon_G,
        ~~n\ge \bar{N},
    \end{align*}
    where $\ca_{\zeta,n}:=({\alpha_n}/{\sqrt{1-2\alpha_{\bar{N}}\varepsilon_\zeta}})
    [{\alpha_nC_{\zeta,1}}/{\beta_n}
    +c_n^2C_{\zeta,3}
    +{\beta_{p,n}^{1/2}C_{\zeta,2}}/{(c_nh_n)}
    +h_n^2C_{\zeta,4}]$
    for some constants $C_{\zeta,1},\ldots,\allowbreak C_{\zeta,4}>0$, $\varepsilon_\zeta:=\varepsilon_\rho-\sqrt{\varepsilon_G}/\delta$.
    Hence, we complete the proof by applying Lemma~\ref{LM_contraction_mapping} and essentially the same argument used in the proof of Proposition~\ref{PROP_rate_measure_uni}.
\endproof


\begin{thebibliography}{42}
\providecommand{\natexlab}[1]{#1}
\providecommand{\url}[1]{\texttt{#1}}
\providecommand{\urlprefix}{URL }

\bibitem[{Acharya et~al.(2017)Acharya, Pedersen, Philippon, \protect\BIBand{}
  Richardson}]{acharya2017measuring}
Acharya VV, Pedersen LH, Philippon T, Richardson M (2017) Measuring systemic
  risk. \emph{The Review of Financial Studies} 30(1):2--47.

\bibitem[{Adrian \protect\BIBand{} Brunnermeier(2008)}]{Adrian2008}
Adrian T, Brunnermeier MK (2008) {CoVaR}. \emph{Federal Reserve Bank of New
  York Staff Report} \urlprefix\url{http://dx.doi.org/10.2139/ssrn.1269446}.

\bibitem[{Ban \protect\BIBand{} Rudin(2019)}]{ban2019big}
Ban GY, Rudin C (2019) The big data newsvendor: Practical insights from machine
  learning. \emph{Operations Research} 67(1):90--108.

\bibitem[{Bertsimas \protect\BIBand{} Kallus(2020)}]{bertsimas2020predictive}
Bertsimas D, Kallus N (2020) From predictive to prescriptive analytics.
  \emph{Management Science} 66(3):1025--1044.

\bibitem[{Bhatnagar et~al.(2001)Bhatnagar, Fu, Marcus, \protect\BIBand{}
  Bhatnagar}]{bhatnagar2001two}
Bhatnagar S, Fu MC, Marcus SI, Bhatnagar S (2001) {Two-timescale algorithms for
  simulation optimization of hidden Markov models}. \emph{IIE Transactions}
  33(3):245--258.

\bibitem[{Bhatnagar et~al.(2003)Bhatnagar, Fu, Marcus, \protect\BIBand{}
  Wang}]{bhatnagar2003two}
Bhatnagar S, Fu MC, Marcus SI, Wang IJ (2003) Two-timescale simultaneous
  perturbation stochastic approximation using deterministic perturbation
  sequences. \emph{ACM Transactions on Modeling and Computer Simulation
  (TOMACS)} 13(2):180--209.

\bibitem[{Bianchi et~al.(2023)Bianchi, De~Luca, \protect\BIBand{}
  Rivieccio}]{bianchi2023non}
Bianchi ML, De~Luca G, Rivieccio G (2023) {Non-Gaussian models for CoVaR
  estimation}. \emph{International Journal of Forecasting} 39(1):391--404.

\bibitem[{Borkar(1997)}]{borkar1997stochastic}
Borkar VS (1997) Stochastic approximation with two time scales. \emph{Systems
  \& Control Letters} 29(5):291--294.

\bibitem[{Borkar(2009)}]{borkar2009stochastic}
Borkar VS (2009) \emph{{Stochastic approximation: A dynamical systems
  viewpoint}}, volume~48 (Springer).

\bibitem[{Cao et~al.(2023)Cao, Hu, \protect\BIBand{} Hu}]{CHH2023}
Cao H, Hu JQ, Hu J (2023) {Black-box CoVaR and its gradient estimation}.
  \emph{Available at SSRN 4583631} .

\bibitem[{Cao et~al.(2025)Cao, Hu, Lian, \protect\BIBand{}
  Yang}]{cao2025infinitesimal}
Cao H, Hu JQ, Lian T, Yang X (2025) {Infinitesimal perturbation analysis (IPA)
  derivative estimation with unknown parameters}. \emph{Automatica} 174:112140.

\bibitem[{Doan(2021)}]{doan2021finite}
Doan TT (2021) Finite-time analysis and restarting scheme for linear
  two-time-scale stochastic approximation. \emph{SIAM Journal on Control and
  Optimization} 59(4):2798--2819.

\bibitem[{Doan(2022)}]{doan2022nonlinear}
Doan TT (2022) {Nonlinear two-time-scale stochastic approximation: Convergence
  and finite-time performance}. \emph{IEEE Transactions on Automatic Control}
  68(8):4695--4705.

\bibitem[{Doan(2025)}]{doan2025fast}
Doan TT (2025) {Fast nonlinear two-time-scale stochastic approximation:
  Achieving ${\cal O}(1/k)$ finite-sample complexity}. \emph{IEEE Transactions
  on Automatic Control} .

\bibitem[{Fan \protect\BIBand{} Hu(1992)}]{fan1992bias}
Fan J, Hu TC (1992) Bias correction and higher order kernel functions.
  \emph{Statistics \& Probability Letters} 13(3):235--243.

\bibitem[{Fan \protect\BIBand{} Yao(1998)}]{fan1998efficient}
Fan J, Yao Q (1998) Efficient estimation of conditional variance functions in
  stochastic regression. \emph{Biometrika} 85(3):645--660.

\bibitem[{Fan \protect\BIBand{} Liu(2016)}]{fan2016direct}
Fan Y, Liu R (2016) A direct approach to inference in nonparametric and
  semiparametric quantile models. \emph{Journal of Econometrics}
  191(1):196--216.

\bibitem[{Glasserman(2004)}]{glasserman2004monte}
Glasserman P (2004) \emph{Monte Carlo methods in financial engineering},
  volume~53 (Springer).

\bibitem[{Haque et~al.(2023)Haque, Khodadadian, \protect\BIBand{}
  Maguluri}]{haque2023tight}
Haque SU, Khodadadian S, Maguluri ST (2023) Tight finite time bounds of
  two-time-scale linear stochastic approximation with markovian noise.
  \emph{arXiv preprint arXiv:2401.00364} .

\bibitem[{H{\"a}rdle \protect\BIBand{}
  M{\"u}ller(1997)}]{hardle1997multivariate}
H{\"a}rdle W, M{\"u}ller M (1997) Multivariate and semiparametric kernel
  regression. Technical report, SFB 373 Discussion Paper.

\bibitem[{Hu \protect\BIBand{} Fu(2025)}]{hu2024}
Hu J, Fu MC (2025) {Technical Note—On the convergence rate of stochastic
  approximation for gradient-based stochastic optimization}. \emph{Operations
  Research} 73(2):1143--1150.

\bibitem[{Hu et~al.(2022)Hu, Peng, Zhang, \protect\BIBand{} Zhang}]{Hu2022}
Hu J, Peng Y, Zhang G, Zhang Q (2022) A stochastic approximation method for
  simulation-based quantile optimization. \emph{INFORMS Journal on Computing}
  34(6):2889--2907, \urlprefix\url{http://dx.doi.org/10.1287/ijoc.2022.1214}.

\bibitem[{Hu et~al.(2025)Hu, Song, \protect\BIBand{} Fu}]{Hu2023}
Hu J, Song M, Fu MC (2025) Quantile optimization via multiple-timescale local
  search for black-box functions. \emph{Operations Research} 73(3):1535--1557.

\bibitem[{Huang et~al.(2024)Huang, Lin, \protect\BIBand{}
  Hong}]{huang2022montecarlo}
Huang W, Lin N, Hong LJ (2024) {Monte Carlo estimation of CoVaR}.
  \emph{Operations Research} 72(6):2337--2357.

\bibitem[{Kallus \protect\BIBand{} Mao(2023)}]{kallus2023stochastic}
Kallus N, Mao X (2023) Stochastic optimization forests. \emph{Management
  Science} 69(4):1975--1994.

\bibitem[{Karimalis \protect\BIBand{} Nomikos(2018)}]{karimalis2018measuring}
Karimalis EN, Nomikos NK (2018) {Measuring systemic risk in the European
  banking sector: A copula CoVaR approach}. \emph{The European Journal of
  Finance} 24(11):944--975.

\bibitem[{Konda \protect\BIBand{} Tsitsiklis(2004)}]{konda2004convergence}
Konda VR, Tsitsiklis JN (2004) Convergence rate of linear two-time-scale
  stochastic approximation. \emph{The Annals of Applied Probability}
  14(2):796--819.

\bibitem[{Kushner \protect\BIBand{} Yin(2003)}]{kushner2003sa}
Kushner HJ, Yin GG (2003) \emph{{Stochastic approximation and recursive
  algorithms and applications}} (Springer), ISBN 0387008942.

\bibitem[{Li \protect\BIBand{} Peng(2024)}]{li2024eliminating}
Li Z, Peng Y (2024) Eliminating ratio bias for gradient-based simulated
  parameter estimation. \emph{arXiv preprint arXiv:2411.12995} .

\bibitem[{Mokkadem \protect\BIBand{} Pelletier(2006)}]{mokkadem2006convergence}
Mokkadem A, Pelletier M (2006) Convergence rate and averaging of nonlinear
  two-time-scale stochastic approximation algorithms. \emph{The Annals of
  Applied Probability} 16(3):1671--1702.

\bibitem[{Parzen(1962)}]{parzen1962estimation}
Parzen E (1962) On estimation of a probability density function and mode.
  \emph{The Annals of Mathematical Statistics} 33(3):1065--1076.

\bibitem[{Qi et~al.(2022)Qi, Cao, \protect\BIBand{}
  Shen}]{qi2022distributionally}
Qi M, Cao Y, Shen ZJ (2022) Distributionally robust conditional quantile
  prediction with fixed design. \emph{Management Science} 68(3):1639--1658.

\bibitem[{Sadana et~al.(2025)Sadana, Chenreddy, Delage, Forel, Frejinger,
  \protect\BIBand{} Vidal}]{sadana2024survey}
Sadana U, Chenreddy A, Delage E, Forel A, Frejinger E, Vidal T (2025) A survey
  of contextual optimization methods for decision-making under uncertainty.
  \emph{European Journal of Operational Research} 320(2):271--289.

\bibitem[{Schucany \protect\BIBand{} Sommers(1977)}]{schucany1977improvement}
Schucany W, Sommers JP (1977) Improvement of kernel type density estimators.
  \emph{Journal of the American Statistical Association} 72(358):420--423.

\bibitem[{Silverman(1978)}]{silverman1978weak}
Silverman BW (1978) Weak and strong uniform consistency of the kernel estimate
  of a density and its derivatives. \emph{The Annals of Statistics}
  6(1):177--184.

\bibitem[{Spall(1992)}]{spall1992multivariate}
Spall JC (1992) Multivariate stochastic approximation using a simultaneous
  perturbation gradient approximation. \emph{IEEE Transactions on Automatic
  Control} 37(3):332--341.

\bibitem[{Tadic \protect\BIBand{} Meyn(2003)}]{tadic2003asymptotic}
Tadic VB, Meyn SP (2003) Asymptotic properties of two time-scale stochastic
  approximation algorithms with constant step sizes. \emph{Proceedings of the
  2003 American Control Conference, 2003.}, volume~5, 4426--4431 (IEEE).

\bibitem[{Wang \protect\BIBand{} Zhang(2025)}]{wang2025statistical}
Wang Y, Zhang X (2025) Statistical inference for weighted sample average
  approximation in contextual stochastic optimization. \emph{arXiv preprint
  arXiv:2503.12747} .

\bibitem[{Wu et~al.(2020)Wu, Zhang, Xu, \protect\BIBand{} Gu}]{wu2020finite}
Wu YF, Zhang W, Xu P, Gu Q (2020) A finite-time analysis of two time-scale
  actor-critic methods. \emph{Advances in Neural Information Processing
  Systems} 33:17617--17628.

\bibitem[{Xu \protect\BIBand{} Cohen(2018)}]{xu2018stock}
Xu Y, Cohen SB (2018) Stock movement prediction from tweets and historical
  prices. \emph{Proceedings of the 56th Annual Meeting of the Association for
  Computational Linguistics}, 1970--1979.

\bibitem[{Yin \protect\BIBand{} Zhang(2005)}]{yin2005discrete}
Yin GG, Zhang Q (2005) \emph{{Discrete-time Markov chains: Two-time-scale
  methods and applications}} (Springer).

\bibitem[{Zhang et~al.(2024)Zhang, Yang, \protect\BIBand{}
  Gao}]{zhang2024optimal}
Zhang L, Yang J, Gao R (2024) Optimal robust policy for feature-based
  newsvendor. \emph{Management Science} 70(4):2315--2329.

\end{thebibliography}

\end{document}